\documentclass[a4paper]{article}
\usepackage{theorem,amsmath,amssymb,amscd}
\usepackage[small]{titlesec}

\theoremstyle{change} \allowdisplaybreaks 

\newcommand{\A}{{\mathbb A}}
\newcommand{\Q}{{\mathbb Q}}
\newcommand{\Z}{{\mathbb Z}}
\newcommand{\R}{{\mathbb R}}
\newcommand{\C}{{\mathbb C}}

\newcommand{\AI}{\mathcal{AI}}

\newcommand{\p}{\mathfrak p}
\newcommand{\IM}{{\rm Im}}

\newcommand{\SH}{\mathfrak h}
\newcommand{\tr}{{\rm tr}}

\newcommand{\GL}{{\rm GL}}
\newcommand{\GU}{\mathrm{GU}}
\newcommand{\PGL}{{\rm PGL}}
\newcommand{\SL}{{\rm SL}}
\newcommand{\SO}{{\rm SO}}
\newcommand{\U}{{\rm U}}
\newcommand{\SSp}{{\rm Sp}}
\newcommand{\GSp}{{\rm GSp}}

\newcommand{\mat}[4]{{\setlength{\arraycolsep}{0.5mm}\left[
 \begin{array}{cc}#1&#2\\#3&#4\end{array}\right]}}

\newcommand{\forget}[1]{}

\newcommand{\qed}{\hspace*{\fill}\rule{1ex}{1ex}}
\newcommand{\nl}{

\vspace{2ex}}
\newcommand{\nll}{

\vspace{1ex}}

\newtheorem{lemma}{Lemma.}[section]
\newtheorem{theorem}[lemma]{Theorem.}
\newtheorem{corollary}[lemma]{Corollary.}
\newtheorem{proposition}[lemma]{Proposition.}
\newtheorem{definition}[lemma]{Definition.}

\begin{document}
\thispagestyle{empty}
\begin{center}
{\bf\Large Bessel models for lowest weight representations of $\GSp(4,\R)$}

\vspace{2ex}
Ameya Pitale\footnote{Department of Mathematics, University of Oklahoma,
Norman, OK 73019-0315, {\tt ameya@math.ou.edu}},
Ralf Schmidt\footnote{Department of Mathematics, University of Oklahoma,
Norman, OK 73019-0315, {\tt rschmidt@math.ou.edu}}

\vspace{5ex}
\begin{minipage}{80ex}
 \small{\sc Abstract.} We prove uniqueness and give precise criteria for existence
 of split and non-split Bessel models for a class of lowest and highest weight
 representations of the groups $\GSp(4,\R)$ and $\SSp(4,\R)$ including all
 holomorphic and anti-holomorphic discrete series representations. Explicit
 formulas for the resulting Bessel functions are obtained by solving systems
 of differential equations. The formulas are applied to derive an integral
 representation for a global $L$-function on $\GSp(4)\times\GL(2)$ involving
 a vector-valued Siegel modular form of degree $2$.
\end{minipage}
\end{center}
\section{Introduction}
Whittaker models for generic representations of a reductive
algebraic group over a local or global field are a very important
tool in representation theory. For non-generic representations of
a classical group over an archimedean or non-archimedean local
field, Bessel models can sometimes provide a substitute for the
missing Whittaker models. Moreover, global Bessel models have been
successfully employed to study $L$-functions and other global
objects, as in \cite{Fu} or \cite{PS}. In view of these
applications, it is desirable to have as much information as
possible about uniqueness and existence of local Bessel models.
\nl In the present paper we prove uniqueness and give precise
criteria for existence of Bessel models for the lowest and highest
weight representations of $\GSp(4,\R)$ and $\SSp(4,\R)$. Our
method is elementary and is based on solving a system of linear
first-order PDEs. This leads not only to uniqueness and existence
criteria, but to explicit formulas for the functions in the Bessel
model. \nl To be a bit more specific, a Bessel model for a
representation of $\GSp(4,\R)$ consists of a space of functions
$B:\:\GSp(4,\R)\rightarrow\C$ that transform on the left in a
certain way under a character of the \emph{Bessel subgroup}, and
such that $\GSp(4,\R)$ acts on this space by right translation.
The Bessel subgroup $R(\R)$ is contained in the Siegel parabolic
subgroup, and can be written as a semidirect product
$R(\R)=T(\R)_0\ltimes U(\R)$, where $U(\R)$ is the unipotent
radical of the Siegel parabolic subgroup, and where $T(\R)$ is the
identity component of the multiplicative group of a quadratic
extension of $\R$. This quadratic extension may be $\R\times\R$,
in which case we speak of a \emph{split} Bessel model; otherwise
we have a \emph{non-split} Bessel model. \nl The class of
representations we consider contains all the holomorphic and
anti-holomorphic discrete series representations, but also the
limits of discrete series representations and certain non-tempered
lowest and highest weight modules. This is exactly the class of
representations that appear in the automorphic representations
generated by holomorphic (scalar or vector valued) Siegel modular
forms. In the split case, the existence question has a simple
answer: None of the lowest weight representations we consider has
a split Bessel model (Corollary
\ref{holdiscserminKBesselsplitthmcor}). The non-split case is more
interesting. In this case $T(\R)\cong\C^\times$, and a Bessel
model exists if and only if (an obvious compatibility with the
central character is satisfied and) the given character $\Lambda$
of $T(\R)$, restricted to the unit circle, is indexed by an
integer $m$ whose absolute value is less than the dimension of the
minimal $K$-type of the lowest weight representation (Theorem
\ref{besselmodeltheorem}). In particular, so-called \emph{special}
Bessel models, i.e., the models for which $\Lambda$ is trivial,
always exist as long as the central character condition is
satisfied. \nl Bessel models for $\GSp(4)$ have also been treated
in the recent preprint \cite{PT}, both in the real and $p$-adic
setting. The overlap with our results consists in the uniqueness
and existence criterion for non-split Bessel models in the
holomorphic discrete series case (see Theorem
\ref{besselmodeltheorem} and \cite{PT}, Theorem 9), which is
obtained in \cite{PT} by completely different methods. Note that
\cite{PT} works with Bessel functionals on the smooth vectors of a
representation, while we work in the category of
$(\mathfrak{g},K)$-modules. \nl Our focus will be on explicit
formulas for Bessel functions, an application of which is given in
Sect.\ \ref{applicationsec}. In Theorem \ref{main-global-thm} we
will obtain an integral representation for the degree-$8$
$L$-function $L(s,\pi_{\mathbf F}\times\tau_f)$ of
$\GSp(4)\times\GL(2)$, where $\mathbf F$ is a vector-valued
holomorphic Siegel modular cusp form of degree $2$ (with respect
to the full modular group) and $f$ is an elliptic Maa{\ss} cusp
form (of arbitrary level). The integral representation is based on
Furusawa's method \cite{Fu} and its extension in \cite{P-S1},
\cite{P-S2}, which involves a global Bessel model for the
automorphic representation $\pi_{\mathbf F}$ generated by $\mathbf
F$. The new aspect here is that $\mathbf F$ is allowed to be
vector-valued, which is possible because of our explicit formulas
for the archimedean Bessel functions. \nl In the first part of
this paper we will define the groups and Lie algebras involved and
introduce the class of lowest and highest weight representations
to be considered. Sections \ref{besselsubgroupssec} and
\ref{besselmodelssec} will make the notion of Bessel model precise
and state some general facts. We start the study of non-split
Bessel models in Sect.\ \ref{nonsplitdoublecosetsec} by recalling
the important double coset decomposition (\ref{HRBesseldecompeq}),
which already appeared in \cite{Fu}. In the next section we derive
formulas for the action of elements of the complexified Lie
algebra on the functions in a Bessel model. These formulas
translate into a system of linear first-order PDEs for the lowest
weight vector $B_0$ in a Bessel model. We solve this system
whenever possible, leading to our first main result Theorem
\ref{holdiscserminKBesselthm} on uniqueness and existence of
certain Bessel functions. In Sect.\ \ref{nonsplitmodelssec} we
translate this result into uniqueness and existence statements for
Bessel models. This requires some more work, but leads to
additional insights into lowest weight modules and their Bessel
models. \nl Split Bessel models will be studied analogously. This
time everything is based on the double coset decomposition
(\ref{disjoint double cosets split}) in Sect.\
\ref{splitdoublecosetsec}. Following that we obtain in Sect.\
\ref{diffopsecsplit} formulas for the action of the complexified
Lie algebra on Bessel functions. We remark that these formulas,
just as their analogues in the non-split case, are independent of
the type of representation considered and are potentially useful
for the study of representations other than the lowest weight
modules considered here. In Sect.\ \ref{hdsrsplitsec} we solve the
resulting system of differential equations in the split case. It
turns out that none of the solutions is of moderate growth,
leading immediately to the non-existence of split Bessel models
for lowest weight representations. \nl Finally, in Sections
\ref{finding-good-vector-sec} -- \ref{intrepsec}, we demonstrate
the applicability of explicit Bessel models by deriving the
integral representation for the $L$-function $L(s,\pi_{\mathbf
F}\times\tau_f)$ mentioned above. The evaluation of the relevant
$p$-adic zeta integrals has been carried out in \cite{Fu},
\cite{P-S1} and \cite{P-S2}, so that we need only evaluate the
archimedean integral.
\section{Definitions and preliminaries}
\subsection{Groups and Lie algebras}\label{grpsec}
Let
$$
 \GSp(4,\R)=\{g\in\GL(4,\R):\:^tgJg=\mu_2(g)J\},\qquad
 J=\begin{bmatrix}&&1\\&&&1\\-1\\&-1\end{bmatrix}.
$$
The function $\mu_2:\:\GSp(4,\R)\rightarrow\R^\times$ is the multiplier homomorphism.
Let $\SSp(4,\R)=\{g\in\GSp(4,\R):\:\mu_2(g)=1\}$.
Let $\mathfrak{g}$ be the Lie algebra of $\GSp(4,\R)$,
and let $\mathfrak{g}^1$ be the Lie algebra of $\SSp(4,\R)$. Then
$\mathfrak{g}=\R\oplus\mathfrak{g}^1$. Explicitly,
$\mathfrak{g}^1=\{\mat{A}{B}{C}{D}\in M(4,\R):\;A=-\,^tD,\;B=\,^tB,\;C=\,^tC\}$.
A basis of $\mathfrak{g}^1$ is given by
\begin{alignat*}{2}
 H_1=&\begin{bmatrix}1&0&0&0\\0&0&0&0\\0&0&-1&0\\0&0&0&0\end{bmatrix}
  ,\qquad\qquad&H_2=&\begin{bmatrix}0&0&0&0\\0&1&0&0\\0&0&0&0
  \\0&0&0&-1\end{bmatrix},\\
 F=&\begin{bmatrix}0&0&1&0\\0&0&0&0\\0&0&0&0\\0&0&0&0\end{bmatrix}
  ,\qquad\qquad&G=&\begin{bmatrix}0&0&0&0\\0&0&0&0\\1&0&0&0
  \\0&0&0&0\end{bmatrix},\\
 R=&\begin{bmatrix}0&0&0&0\\0&0&0&1\\0&0&0&0\\0&0&0&0\end{bmatrix}
  ,\qquad\qquad&R'=&\begin{bmatrix}0&0&0&0\\0&0&0&0\\0&0&0&0
  \\0&1&0&0\end{bmatrix},\\
 P=&\begin{bmatrix}0&0&0&0\\1&0&0&0\\0&0&0&-1\\0&0&0&0\end{bmatrix}
  ,\qquad\qquad&P'=&\begin{bmatrix}0&1&0&0\\0&0&0&0\\0&0&0&0
  \\0&0&-1&0\end{bmatrix},\\
 Q=&\begin{bmatrix}0&0&0&1\\0&0&1&0\\0&0&0&0\\0&0&0&0\end{bmatrix}
  ,\qquad\qquad&Q'=&\begin{bmatrix}0&0&0&0\\0&0&0&0\\0&1&0&0
  \\1&0&0&0\end{bmatrix}.
\end{alignat*}
\nll
A convenient basis for the complexified Lie algebra $\mathfrak{g}^1_\C=
\mathfrak{g}^1\otimes\C$ is as follows.
\begin{alignat*}{2}
 &Z=-i\begin{bmatrix}0&0&1&0\\0&0&0&0\\-1&0&0&0\\0&0&0&0\end{bmatrix}
  ,\qquad\qquad&&Z'=-i\begin{bmatrix}0&0&0&0\\0&0&0&1\\0&0&0&0
  \\0&-1&0&0\end{bmatrix},\\
 &N_+=\frac12\begin{bmatrix}0&1&0&-i\\-1&0&-i&0\\
  0&i&0&1\\i&0&-1&0\end{bmatrix}
  ,\qquad\qquad&&N_-=\frac12\begin{bmatrix}
  0&1&0&i\\-1&0&i&0\\0&-i&0&1\\-i&0&-1&0\end{bmatrix},\\
 &X_+=\frac12\begin{bmatrix}1&0&i&0\\0&0&0&0\\i&0&-1&0\\
  0&0&0&0\end{bmatrix}
  ,\qquad\qquad&&X_-=\frac12\begin{bmatrix}
  1&0&-i&0\\0&0&0&0\\-i&0&-1&0\\0&0&0&0\end{bmatrix},\\
 &P_{1+}=\frac12\begin{bmatrix}0&1&0&i\\1&0&i&0\\
  0&i&0&-1\\i&0&-1&0\end{bmatrix}
  ,\qquad\qquad&&P_{1-}=\frac12\begin{bmatrix}
  0&1&0&-i\\1&0&-i&0\\0&-i&0&-1\\-i&0&-1&0\end{bmatrix},\\
 &P_{0+}=\frac12\begin{bmatrix}0&0&0&0\\0&1&0&i\\
  0&0&0&0\\0&i&0&-1\end{bmatrix}
  ,\qquad\qquad&&P_{0-}=\frac12\begin{bmatrix}
  0&0&0&0\\0&1&0&-i\\0&0&0&0\\0&-i&0&-1\end{bmatrix}.
\end{alignat*}
The following multiplication table for this basis will be useful.
$$\renewcommand{\arraystretch}{1.5}
\begin{array}{c|cccccccccc}
 &Z&Z'&N_+&N_-&X_+&X_-&P_{1+}&P_{1-}&P_{0+}&P_{0-}\\\hline
 Z&0&0&N_+&-N_-&2X_+&-2X_-&P_{1+}&-P_{1-}&0&0\\
 Z'&0&0&-N_+&N_-&0&0&P_{1+}&-P_{1-}&2P_{0+}&-2P_{0-}\\
 N_+&-N_+&N_+&0&Z'-Z&0&-P_{1-}&2X_+&-2P_{0-}&P_{1+}&0\\
 N_-&N_-&-N_-&Z-Z'&0&-P_{1+}&0&-2P_{0+}&2X_-&0&P_{1-}\\
 X_+&-2X_+&0&0&P_{1+}&0&Z&0&N_+&0&0\\
 X_-&2X_-&0&P_{1-}&0&-Z&0&N_-&0&0&0\\
 P_{1+}&-P_{1+}&-P_{1+}&-2X_+&2P_{0+}&0&-N_-&0&Z+Z'&0&N_+\\
 P_{1-}&P_{1-}&P_{1-}&2P_{0-}&-2X_-&-N_+&0&-Z-Z'&0&N_-&0\\
 P_{0+}&0&-2P_{0+}&-P_{1+}&0&0&0&0&-N_-&0&Z'\\
 P_{0-}&0&2P_{0-}&0&-P_{1-}&0&0&-N_+&0&-Z'&0\\
\end{array}
$$
Let $\mathfrak{h}_\R$ be the real subspace spanned by $Z$ and $Z'$, and
let $\mathfrak{h}_\C$ be its complexification. Identifying an $\R$-linear map
$\lambda:\:\mathfrak{h}_\R\rightarrow\R$ with the pair $(\lambda(Z),\lambda(Z'))$,
we get an isomorphism $\mathfrak{h}_\R^*\cong\R^2$.
Such a map $\lambda$ is \emph{analytically integral} if $(\lambda(Z),\lambda(Z'))\in\Z^2$.
The root system is $\Delta=\{(\pm2,0),(0,\pm2),(\pm1,\pm1),(\pm1,\mp1)\}$.
The following diagram indicates the analytically integral elements, as well as the
roots and the elements of the Lie algebra spanning the corresponding root spaces.
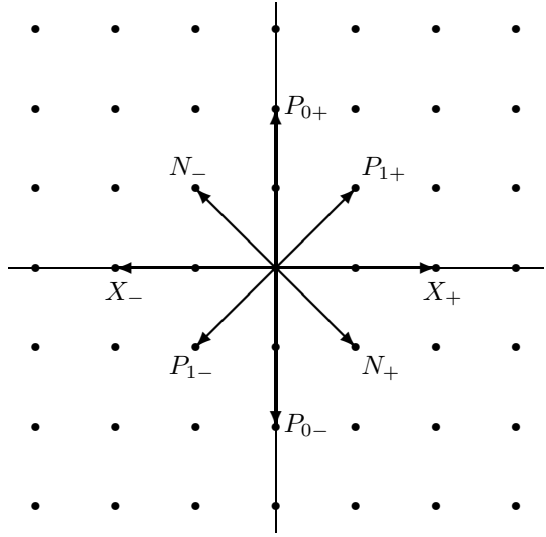
\begin{center}
\begin{picture}(200, 200)(-100, -100)
\put(-100,0){\line(1,0){200}}
\put(0,100){\line(0,-1){200}}

\thicklines
\put(0,0){\vector(1,0){60}}
\put(0,0){\vector(-1,0){60}}
\put(0,0){\vector(0,1){60}}
\put(0,0){\vector(0,-1){60}}
\put(0,0){\vector(1,1){30}}
\put(0,0){\vector(-1,-1){30}}
\put(0,0){\vector(1,-1){30}}
\put(0,0){\vector(-1,1){30}}

\put(55,-12){$X_+$}
\put(-64,-12){$X_-$}
\put(32,-40){$N_+$}
\put(-40,35){$N_-$}
\put(3,58){$P_{0+}$}
\put(3,-62){$P_{0-}$}
\put(32,35){$P_{1+}$}
\put(-40,-40){$P_{1-}$}

\put(-90,90){\circle*{3}}
\put(-60,90){\circle*{3}}
\put(-30,90){\circle*{3}}
\put(0,90){\circle*{3}}
\put(30,90){\circle*{3}}
\put(60,90){\circle*{3}}
\put(90,90){\circle*{3}}

\put(-90,60){\circle*{3}}
\put(-60,60){\circle*{3}}
\put(-30,60){\circle*{3}}
\put(0,60){\circle*{3}}
\put(30,60){\circle*{3}}
\put(60,60){\circle*{3}}
\put(90,60){\circle*{3}}

\put(-90,30){\circle*{3}}
\put(-60,30){\circle*{3}}
\put(-30,30){\circle*{3}}
\put(0,30){\circle*{3}}
\put(30,30){\circle*{3}}
\put(60,30){\circle*{3}}
\put(90,30){\circle*{3}}

\put(-90,0){\circle*{3}}
\put(-60,0){\circle*{3}}
\put(-30,0){\circle*{3}}
\put(0,0){\circle*{3}}
\put(30,0){\circle*{3}}
\put(60,0){\circle*{3}}
\put(90,0){\circle*{3}}

\put(-90,-30){\circle*{3}}
\put(-60,-30){\circle*{3}}
\put(-30,-30){\circle*{3}}
\put(0,-30){\circle*{3}}
\put(30,-30){\circle*{3}}
\put(60,-30){\circle*{3}}
\put(90,-30){\circle*{3}}

\put(-90,-60){\circle*{3}}
\put(-60,-60){\circle*{3}}
\put(-30,-60){\circle*{3}}
\put(0,-60){\circle*{3}}
\put(30,-60){\circle*{3}}
\put(60,-60){\circle*{3}}
\put(90,-60){\circle*{3}}

\put(-90,-90){\circle*{3}}
\put(-60,-90){\circle*{3}}
\put(-30,-90){\circle*{3}}
\put(0,-90){\circle*{3}}
\put(30,-90){\circle*{3}}
\put(60,-90){\circle*{3}}
\put(90,-90){\circle*{3}}

\end{picture}
\end{center}
Let $\mathfrak{k}$ be the $1$-eigenspace of the Cartan involution
$\theta:\:X\mapsto -X^t$, and let $\p$ the $(-1)$-eigenspace of $\theta$. Then
$$
 \mathfrak{k}_\C:=\mathfrak{k}\otimes\C=\langle Z,Z',N_+,N_-\rangle,
$$
and $\p_\C:=\p\otimes\C=\p_++\p_-$, with the maximal abelian subalgebras
$$
 \p_+=\langle X_+,P_{1+},P_{0+}\rangle,\qquad\qquad
 \p_-=\langle X_-,P_{1-},P_{0-}\rangle.
$$
The decomposition $\mathfrak{g}^1=\mathfrak{k}+\p_++\p_-$
holds, and $\mathfrak{k}+\p_\pm$ is a parabolic subalgebra of $\mathfrak{g}^1$.
The compact roots are $\Delta_c=\{(\pm1,\mp1)\}$, and the non-compact roots are
$\Delta_n=\{(\pm2,0),(0,\pm2),(\pm1,\pm1)\}$. The Weyl group $W$ has eight elements
and is generated by the reflections about the hyperplanes orthogonal to the root vectors.
The compact Weyl group $W_K$ has two elements and is generated by the reflection
about the hyperplane orthogonal to the compact roots.
\subsection{The maximal compact subgroup}
Let $K$ be the standard maximal compact subgroup of $\GSp(4,\R)$,
and let $K^1$ be the standard maximal compact subgroup of $\SSp(4,\R)$.
Then $K^1$ is the identity component of $K$, and has index $2$ in $K$. If
$\SH_2=\big\{Z\in M(2,\C):\;Z\mbox{ symmetric and }\IM(Z)
\mbox{ positive definite}\big\}$
is the Siegel upper half plane of degree $2$, and if
$\SSp(4,\R)$ acts on $\SH_2$ by fractional linear transformations,
$$
 \mat{A}{B}{C}{D}\langle Z \rangle := (AZ+B)(CZ+D)^{-1},
$$
then $K^1$ is the stabilizer of the element $I=\mat{i}{}{}{i}\in\SH_2$.
It is easy to check that $K^1\simeq\U(2)$ via
$$
 K^1\ni\mat{A}{B}{-B}{A}\longmapsto A+iB\in\U(2).
$$
Let $J(h,Z) := CZ+D$. If we let $\mathfrak k\subset\mathfrak g$ be the Lie algebra
of $K^1$, then $\mathfrak k$ coincides with the Lie algebra of the same name
mentioned in the previous section.
\subsubsection*{Representations of $K^1$}
We say a vector $v$ in a representation of $\mathfrak{k}_\C$
has \emph{weight} $(l,l')\in\Z^2$ if $Zv=lv$ and $Z'v=l'v$.
To describe the irreducible representations of $K$, we first describe
those of $K^1$. The latter are in one-one correspondence with the
irreducible representations of $\mathfrak k_\C$ with integral weights.
This Lie algebra is a direct sum
$$
 \mathfrak k_\C=\langle Z-Z',N_+,N_-\rangle\oplus\langle Z+Z'\rangle,
$$
where $Z+Z'$ spans the center of $\mathfrak k_\C$. Now the subalgebra
$\langle Z-Z',N_+,N_-\rangle$ is isomorphic to $\mathfrak{su}(2)$,
and therefore its irreducible representations are indexed by non-negative
integers indicating the weight of a highest weight vector (a vector
annihilated by $N_+$). We also can prescribe any integer with which
$Z+Z'$ is supposed to act, but have to make sure the resulting representation
of $\mathfrak k_\C$ has integral weights. This shows that
the (isomorphism classes of) irreducible representations of $K^1$
are in one-one correspondence with the set $\{(l,l')\in\Z^2:\;l\geq l'\}$, or,
in other words, with the analytically integral elements of $\mathfrak{h}_\R^*$
modulo the action of $W_K$.
If we let $\rho_{l,l'}$ be the representation corresponding to the
pair $(l,l')$, then $\rho_{l,l'}$ is characterized by the property
that it contains a non-zero vector of weight $(l,l')$ that is
annihilated by $N_+$.
The weight structure of $\rho_{l,l'}$ is symmetric with respect to the main diagonal
(the wall orthogonal to the compact roots).
It contains a highest weight vector of weight $(l,l')$ (annihilated
by $N_+$), and a lowest weight vector of weight $(l',l)$ (annihilated
by $N_-$). It contains the weights ``between'' these two extremes
exactly once. The one-dimensional representations are the $\rho_{l,l}$
with $l\in\Z$. The representation $\rho_{0,0}$ is the trivial one.
Evidently, the dimension of $\rho_{l,l'}$ is $l-l'+1$.
The $\rho_{l,l'}$ with trivial central character are those for which $l,l'$ are
both even or both odd. These representations are odd-dimensional.
\subsubsection*{Representations of $K$}
It is now easy to describe the representations of $K$. They are
all obtained by induction from representations of $K^1$. The induction
process has the effect of making the weight structure point symmetric
with respect to the origin. More precisely, if $(\rho,V)$ is a representation
of $K$, and if $v\in V$ has weight $(l,l')$, then
$\rho({\rm diag}(1,1,-1,-1))v$ has weight $(-l,-l')$. Thus the
weight structure of an irreducible representation of $K$ combines
that of $\rho_{l,l'}$ and $\rho_{-l',-l}$, for some pair $(l,l')$.
The representations $\rho_{l,-l}$ of $K^1$ with $l\geq0$ extend in two different ways
to representations of $K$.
\subsubsection*{Coordinates on $K^1$}
The following coordinates on $K^1$ will be convenient. We map
\begin{equation}\label{Kcoordeq}
 \R^4\ni(\varphi_1,\varphi_2,\varphi_3,\varphi_4)\longmapsto
 r_1(\varphi_1)r_2(\varphi_2)r_3(\varphi_3)r_4(\varphi_4)\in K^1,
\end{equation}
where
\begin{alignat*}{2}
 &r_1(\varphi_1)=\begin{bmatrix}\cos(\varphi_1)&\sin(\varphi_1)\\
  -\sin(\varphi_1)&\cos(\varphi_1)\\&&\cos(\varphi_1)&\sin(\varphi_1)\\
  &&-\sin(\varphi_1)&\cos(\varphi_1)\end{bmatrix},&\qquad
 &r_3(\varphi_3)=\begin{bmatrix}\cos(\varphi_3)&&\sin(\varphi_3)&\\
  &1&&\\-\sin(\varphi_3)&&\cos(\varphi_3)&\\&&&1\end{bmatrix},\\
 &r_2(\varphi_2)=\begin{bmatrix}\cos(\varphi_2)&&&\sin(\varphi_2)\\
  &\cos(\varphi_2)&\sin(\varphi_2)&\\&-\sin(\varphi_2)&\cos(\varphi_2)&\\
  -\sin(\varphi_2)&&&\cos(\varphi_2)\end{bmatrix},&
 &r_4(\varphi_4)=\begin{bmatrix}1&&&\\&\cos(\varphi_4)&&\sin(\varphi_4)\\
  &&1&\\&-\sin(\varphi_4)&&\cos(\varphi_4)\end{bmatrix}.
\end{alignat*}
One can check that the differential of this map at $(0,0,0,0)$
is regular, so that we get coordinates on $K^1$ in a
neighborhood of the identity. This will suffice: We will only be
dealing with $K^1$-finite functions, which are analytic, and therefore
determined on a small neighborhood of the identity of the connected
group $K^1$.
\nl
If $l,l'$ are integers, and a function $\Psi$ on $K^1$ has the property
$$
 \Psi(gr_3(\varphi_3)r_4(\varphi_4))=e^{i(l\varphi_3+l'\varphi_4)}
   \Psi(g)\qquad\text{for all }\varphi_3,\varphi_4\in\R,
$$
then we say that $\Psi$ has \emph{weight} $(l,l')$. This is
equivalent to our previous notion of weight, if $\Psi$ is an element
of a function space (e.g.\ $L^2(K^1)$) on which $K^1$ acts by right
translation.
\subsection{Lowest and highest weight representations}\label{lowestweightsec}
We will describe certain lowest and highest weight representations of the
group $\SSp(4,\R)$. Let $l,l'$ be integers
with $l\geq l'>0$. Then there exists an irreducible representation
$\mathcal{E}(l,l')$ of $\SSp(4,\R)$ characterized by the existence of a non-zero vector $v$
of weight $(l,l')$ with the property
\begin{equation}\label{typeIannpropeq}
 N_+v=X_-v=P_{1-}v=P_{0-}v=0.
\end{equation}
If $l'\geq3$, then $\mathcal{E}(l,l')$ is a holomorphic discrete series representation
with Harish-Chandra parameter $\lambda=(l-1,l'-2)$.
If $l'=2$, then $\mathcal{E}(l,l')$ is a limit of discrete
series representations. If $l'=1$, then $\mathcal{E}(l,l')$ is not a discrete
series or limit of discrete series representation.
\nl
Similarly, if $l'\leq l<0$, then there exists an irreducible representation
$\mathcal{E}(l,l')$ characterized by the existence of a non-zero vector $v$
of weight $(l,l')$ with the property
\begin{equation}\label{typeIVannpropeq}
 N_+v=X_+v=P_{1+}v=P_{0+}v=0.
\end{equation}
If $l\leq-3$, then $\mathcal{E}(l,l')$ coincides with the anti-holomorphic discrete series
representation with Harish-Chandra parameter $\lambda=(l+2,l'+1)$.
If $l=-2$, then $\mathcal{E}(l,l')$ is a limit of discrete series representations.
If $l=-1$, then $\mathcal{E}(l,l')$ is not a discrete
series or limit of discrete series representation.
\nl
The following diagrams illustrate the regions containing all the $K^1$-types of
$\mathcal{E}(l,l')$ in the holomorphic (left) resp.\ anti-holomorphic (right) case.
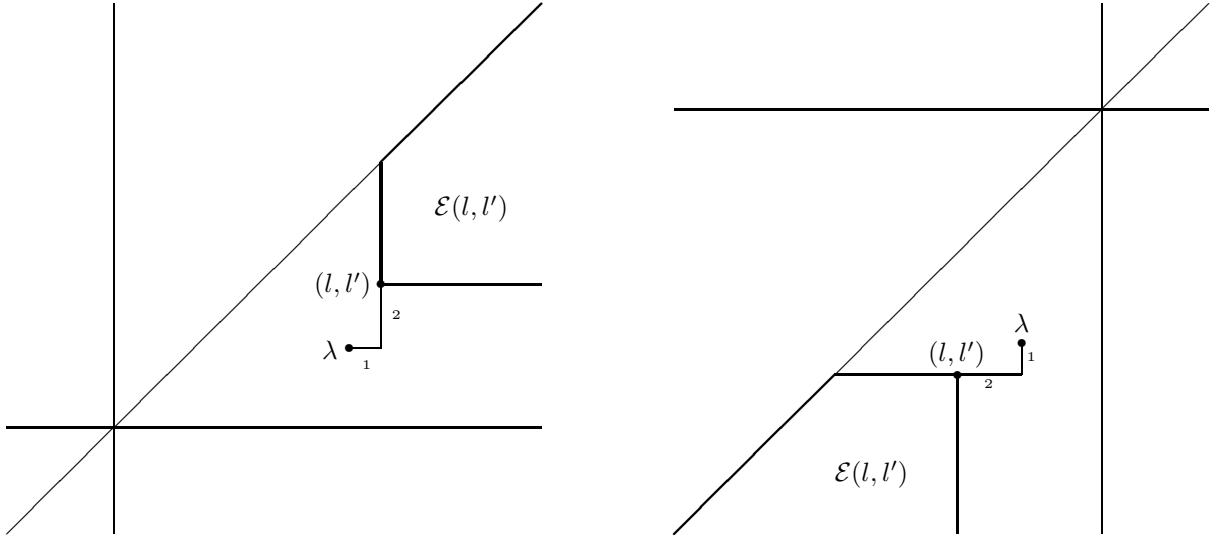
\begin{center}
\begin{picture}(200, 200)(-40, -40)
\put(-40,0){\line(1,0){200}}
\put(0,-40){\line(0,1){200}}
\put(-40,-40){\line(1,1){200}}

\put(100,54){\circle*{3}}
\put(88,30){\circle*{3}}
\put(88,30){\line(1,0){12}}
\put(100,30){\line(0,1){24}}
\put(93,22){\tiny$1$}
\put(104,41){\tiny$2$}
\put(78,26){$\lambda$}
\put(75,51){$(l,l')$}

\thicklines
\put(100,54){\line(1,0){60}}
\put(100,54){\line(0,1){46}}
\put(100,100){\line(1,1){60}}
\put(120,80){$\mathcal{E}(l,l')$}

\end{picture}
\hspace{10ex}
\begin{picture}(200, 200)(-160, -160)
\put(-160,0){\line(1,0){200}}
\put(0,-160){\line(0,1){200}}
\put(-160,-160){\line(1,1){200}}

\put(-54,-100){\circle*{3}}
\put(-30,-88){\circle*{3}}
\put(-30,-100){\line(-1,0){24}}
\put(-30,-100){\line(0,1){12}}
\put(-28,-95){\tiny$1$}
\put(-44,-105){\tiny$2$}
\put(-33,-84){$\lambda$}
\put(-65,-95){$(l,l')$}

\thicklines
\put(-54,-100){\line(0,-1){60}}
\put(-54,-100){\line(-1,0){46}}
\put(-100,-100){\line(-1,-1){60}}
\put(-100,-140){$\mathcal{E}(l,l')$}

\end{picture}
\end{center}
One reason the lowest weight representations $\mathcal{E}(l,l')$ are important
is that they appear as
the archimedean components of the automorphic representations generated by Siegel modular
forms of degree $2$. If $l=l'$, then they correspond to scalar-valued Siegel modular
forms of weight $l$. If $l>l'$, then they correspond to vector-valued Siegel modular
forms; see \cite{AS}. Using our formulas for Bessel models for $\mathcal{E}(l,l')$ obtained
further below, we will demonstrate in Sect.\ \ref{applicationsec} how to
obtain an integral representation for an $L$-function on
$\GSp(4)\times\GL(2)$ involving a vector-valued Siegel modular form.
\subsection{Representations of $\SSp(4,\R)$ and of $\GSp(4,\R)$}\label{SpGSpsec}
Let $\SSp(4,\R)^\pm=\{g\in\GSp(4,\R):\:\mu_2(g)=\pm1\}$. We fix the element
$\epsilon={\rm diag}(1,1,-1,-1)$. Then $\SSp(4,\R)^\pm=\SSp(4,\R)\sqcup
\epsilon \SSp(4,\R)$. We work in the category of $(\mathfrak{g},K)$-modules
for $\GSp(4,\R)$, resp.\ $(\mathfrak{g}^1,K)$-modules for $\SSp(4,\R)^\pm$,
resp.\ $(\mathfrak{g}^1,K^1)$-modules for $\SSp(4,\R)$. If $(\pi,V)$ is a $(\mathfrak{g}^1,K^1)$-module of $\SSp(4,\R)$, then we define
another $(\mathfrak{g}^1,K^1)$-module $(\pi^\epsilon,V^\epsilon)$ by
$V^\epsilon=V$ and
$$
 \pi^\epsilon(X)=\pi({\rm Ad}(\epsilon)X)\qquad(X\in\mathfrak{g}^1),\qquad
 \pi^\epsilon(k)=\pi(\epsilon k\epsilon^{-1})\qquad(k\in K^1).
$$
If $(l,l')$ is a weight for $\pi$, then $(-l,-l')$ is a weight for $\pi^\epsilon$.
\subsubsection*{Representations of $\SSp(4,\R)^\pm$}
Assume that $\pi$ is an irreducible $(\mathfrak{g}^1,K^1)$-module.
If $\pi^\varepsilon\cong\pi$, then $\pi$ can be extended in exactly two
non-isomorphic ways to a $(\mathfrak{g}^1,K)$-module for $\SSp(4,\R)^\pm$. If
$\pi^\epsilon\not\cong\pi$, then we extend the $(\mathfrak{g}^1,K^1)$-module
structure on the direct sum $V\oplus V^\epsilon$ to a $(\mathfrak{g}^1,K)$-module
structure by requiring that $\pi(\epsilon)(v_1,v_2)=(v_2,v_1)$. This
$(\mathfrak{g}^1,K)$-module is irreducible. We denote it by $I(\pi)$.
Every irreducible $(\mathfrak{g}^1,K)$-module is either of the form $I(\pi)$,
or is obtained by extending an $\epsilon$-invariant $(\mathfrak{g}^1,K^1)$-module.
\nl
{\bf Example:} Let $\pi=\mathcal{E}(l,l')$ with $l\geq l'>0$ be a lowest weight
representation of $\SSp(4,\R)$, as above. We think of $\mathcal{E}(l,l')$ as the
underlying $(\mathfrak{g}^1,K^1)$-module. Then $\pi^\epsilon$ is the highest
weight representation $\mathcal{E}(-l',-l)$. The induced $(\mathfrak{g}^1,K)$-module
combines the space of $\mathcal{E}(l,l')$ and the space of
$\mathcal{E}(-l',-l)$. We denote this $(\mathfrak{g}^1,K)$-module by
$\mathcal{E}(l,l')$.
\subsubsection*{Representations of $\GSp(4,\R)$}
Given a complex number $s$, we can extend a representation of $\SSp(4,\R)^\pm$
to a representation of $\GSp(4,\R)^\pm\cong\R_{>0}\times\SSp(4,\R)^\pm$ by requiring that
${\rm diag}(\gamma,\gamma,\gamma,\gamma)$, $\gamma>0$, acts by multiplication with
$\gamma^s$. On the level of Lie algebras, the central element ${\rm diag}(1,1,1,1)\in
\mathfrak{g}$ acts by multiplication with $s$. If $\pi$ is a
$(\mathfrak{g}^1,K^1)$-module with $\pi\not\cong\pi^\epsilon$, and $I(\pi)$ is
the irreducible $(\mathfrak{g}^1,K)$-module constructed from $\pi$, then we denote
by $I_s(\pi)$ the extension to a representation of $\GSp(4,\R)$.
If $\mathcal{E}(l,l')$, $l\geq l'>0$, is
one of the lowest weight $(\mathfrak{g}^1,K)$-modules described in the
previous paragraph, then we denote by $\mathcal{E}_s(l,l')$
the corresponding $(\mathfrak{g},K)$-module. We call these modules
lowest weight representations of $\GSp(4,\R)$ (even though they have both a lowest
and a highest weight, and even though they are not
representations of $\GSp(4,\R)$ at all).
\subsection{Bessel subgroups}\label{besselsubgroupssec}
Let $U(\R)=\{\mat{1}{X}{}{1}\in\GSp(4,\R):\:^tX=X\}$.
Let $S$ be a non-degenerate real symmetric matrix, and let $\theta$ be the character of
$U(\R)$ given by $\theta(\mat{1}{X}{}{1})=e^{2\pi i{\rm tr}(SX)}$. Explicitly, if
$S=\mat{a}{\frac b2}{\frac b2}{c}$, then
\begin{equation}\label{thetadef2eq}
 \theta(\begin{bmatrix}1&&x&y\\&1&y&z\\&&1\\&&&1\end{bmatrix})=e^{2\pi i(ax+by+cz)}.
\end{equation}
Let $L=\{\mat{g}{}{}{\lambda\,^tg^{-1}}:\:g\in\GL(2,\R),
\:\lambda\in\R^\times\}$ be the Levi component of the Siegel parabolic subgroup
of $\GSp(4,\R)$. We always think of $\GL(2,\R)$
embedded as a subgroup of $L$ via $g\mapsto\mat{g}{}{}{\det(g)\,^tg^{-1}}$. Let
\begin{equation}\label{TRthetaeq1}
 T(\R)=\{g\in L:\:\theta(u)=\theta(gug^{-1})\:\text{for all }u\in U(\R)\}
\end{equation}
and $T^1(\R)=T(\R)\cap\SSp(4,\R)$. Then $T(\R)=\{g\in\GL(2,\R):\:^tgSg=\det(g)S\}$
and $T^1(\R)=\{g\in\SL(2,\R):\:^tgSg=S\}$. Let $T(\R)_0$ (resp.\ $T^1(\R)_0$) be
the identity component of $T(\R)$ (resp.\ $T^1(\R)$) in the real topology.
Then $T(\R)_0$ contains $\{\mat{a}{}{}{a}:\:a>0\}\cong\R_{>0}$, which, as a subgroup
of $\GSp(4,\R)$, corresponds to central elements with positive diagonal entries.
We have
\begin{equation}\label{TR0T1R0eq}
 T(\R)_0=\R_{>0}\times T^1(\R)_0.
\end{equation}
We consider two special cases.
\begin{enumerate}
 \item Let $S=\pm\mat{1}{0}{0}{1}$ (the definite case). Then
  \begin{equation}
   T(\R) = \{\mat{x}{y}{-y}{x}:\:x,y\in\R,\:x^2+y^2>0\}\cong\C^\times
  \end{equation}
  via $\mat{x}{y}{-y}{x}\mapsto x+iy$. The subgroup $T^1(\R)$ corresponds to
  elements of the unit circle. In particular, $T(\R)$ and $T^1(\R)$ are connected.
 \item Let $S=\mat{0}{1}{1}{0}$ (the split case). Then
  \begin{equation}
   T(\R) = \{\mat{x}{}{}{y}:\:x,y\in\R^\times\}\cong\R^\times\times\R^\times.
  \end{equation}
  This group has four connected components, with $T(\R)_0\cong\R_{>0}\times\R_{>0}$.
  The group $T^1(\R)\cong\R^\times$ has two connected components, with
  $T^1(\R)_0\cong\R_{>0}$.
\end{enumerate}
For any $S$, let
\begin{equation}\label{besselsubgroupdefeq}
 R(\R)=T(\R)_0U(\R),\qquad R^1(\R)=T^1(\R)_0U(\R).
\end{equation}
We call $R(\R)$ (resp.\ $R^1(\R)$) the \emph{Bessel subgroup} of $\GSp(4,\R)$
(resp.\ $\SSp(4,\R)$) with respect to $S$. Evidently, $R(\R)=\R_{>0}\times R^1(\R)$.
Summarizing, we have the following groups and subgroups.
\begin{equation}\label{groupsdiagrameq}\renewcommand{\arraystretch}{1.2}
\begin{array}{ccc}
 \GSp(4,\R)&=&\R_{>0}\times\SSp(4,\R)^\pm\\
 \uparrow&&\uparrow\\
 \GSp(4,\R)^+&=&\R_{>0}\times\SSp(4,\R)\\
 \uparrow&&\uparrow\\
 R(\R)&=&\R_{>0}\times R^1(\R)\\
 \uparrow&&\uparrow\\
 T(\R)_0&=&\R_{>0}\times T^1(\R)_0
\end{array}
\end{equation}
\subsection{Bessel models}\label{besselmodelssec}
If $G$ is a Lie group, then $G$ acts on the space of smooth functions
$\Phi:\:G\rightarrow\C$ by right translation: $(h.\Phi)(g)=\Phi(gh)$.
The Lie algebra of $G$ acts on the same space via the derived representation,
$$
 (X.\Phi)(g)=\frac d{dt}\Big|_0\Phi(g\exp(tX)).
$$
We call this action of the Lie algebra also \emph{right translation}.
\subsubsection*{Growth condition}
We define a norm function $\|\cdot\|$ on $\GSp(4,\R)$ by
\begin{equation}\label{normdefeq}
 \|g\|=\Big(\mu_2(g)^{-2}+\sum_{i,j=1}^4g_{ij}^2\Big)^{1/2},\qquad
 g=(g_{ij})\in\GSp(4,\R),
\end{equation}
where $\mu_2$ denotes the  multiplier. We say a function
$\Phi:\:\GSp(4,\R)\rightarrow\C$ or $\Phi:\:\SSp(4,\R)\rightarrow\C$
is \emph{slowly increasing}
(or of \emph{moderate growth}) if there exist
positive constants $\alpha,\beta$ such that
\begin{equation}\label{slowlyincreasingeq1}
 |\Phi(g)|\leq\alpha\|g\|^\beta\qquad\text{for all }g\in\GSp(4,\R)\quad
 (\text{resp. }g\in\SSp(4,\R)).
\end{equation}
The norm is designed so that the function $|\mu_2(g)|$ is slowly increasing.
A necessary condition for $\Phi$ to be slowly increasing is that there exist
positive constants $\alpha,\beta$ such that
\begin{equation}\label{slowlyincreasingeq2}
 |B(\begin{bmatrix}\lambda\zeta\\&\lambda\zeta^{-1}\\&&\zeta^{-1}\\&&&\zeta
  \end{bmatrix})|\leq \alpha(\lambda\zeta)^\beta\qquad\text{for all }
  \lambda,\zeta\in\R^\times,\;|\lambda|,\zeta>1.
\end{equation}
Evidently, if $\Phi:\:\GSp(4,\R)\rightarrow\C$ is slowly increasing, then
its restriction to $\SSp(4,\R)$ is also slowly increasing.
\subsubsection*{Definition of Bessel models}
Let $\Lambda$ be a character of $T(\R)_0$.
In the definite case, $\Lambda$ is a character of $\C^\times$, and in the
split case, $\Lambda$ is a character of $\R_{>0}\times\R_{>0}$.
Let $\Lambda^1$ be the restriction of $\Lambda$ to $T^1(\R)_0$.
In view of (\ref{TR0T1R0eq}), every character of $T^1(\R)_0$ is obtained
by such a restriction.
Since the elements $g$ of $T(\R)_0$ satisfy $\theta(gug^{-1})=\theta(u)$ for
all $u\in U(\R)$, the map
$$
 gu\longmapsto\Lambda(g)\theta(u),\qquad g\in T(\R)_0,\;u\in U(\R),
$$
defines a character of the Bessel subgroup $R(\R)$, which we denote by
$\Lambda\otimes\theta$. Its restriction to $R^1(\R)$ is a character denoted
by $\Lambda^1\otimes\theta$.
\nl
Let $\mathcal{S}(\Lambda,\theta)$ be the space of functions
$B:\:\GSp(4,\R)\rightarrow\C$ with the following properties.
\begin{enumerate}
 \item $B$ is smooth and $K$-finite.
 \item $B(tug)=\Lambda(t)\theta(u)B(g)$ for all $t\in T(\R)_0$, $u\in U(\R)$
  and $g\in\GSp(4,\R)$.
 \item $B$ is slowly increasing.
\end{enumerate}
Let $\mathcal{S}^1(\Lambda^1,\theta)$ be the space of functions
$B:\:\SSp(4,\R)\rightarrow\C$ with the following properties.
\begin{enumerate}
 \item $B$ is smooth and $K$-finite.
 \item $B(tug)=\Lambda^1(t)\theta(u)B(g)$ for all $t\in T^1(\R)_0$, $u\in U(\R)$
  and $g\in\SSp(4,\R)$.
 \item $B$ is slowly increasing.
\end{enumerate}
It is clear that restriction defines a linear map
$\mathcal{S}(\Lambda,\theta)\rightarrow\mathcal{S}^1(\Lambda^1,\theta)$.
We claim that this map is onto. Indeed, let $B\in\mathcal{S}^1(\Lambda^1,\theta)$
be given. First extend $B$ to a function on $\SSp(4,\R)^\pm$ by setting it
equal to zero on elements of negative multiplier. This function satisfies the
correct transformation property under elements of $R^1(\R)$.
Since $\GSp(4,R)=\R_{>0}\times\SSp(4,\R)^\pm$ and $R(\R)=\R_{>0}\times R^1(\R)$,
we can extend $B$ further to a function on $\GSp(4,\R)$ satisfying the
correct transformation property under $R(\R)$; see (\ref{groupsdiagrameq}).
\begin{definition}\label{besselmodeldef}
 Let $S$ be a non-degenerate real symmetric matrix, and let $\theta$ be the
 corresponding character of $U(\R)$, as above.
 \begin{enumerate}
  \item Let $(\pi,V)$ be a $(\mathfrak{g},K)$-module.
   Let $\Lambda$ be a character of $T(\R)_0$. A \emph{$(\Lambda,\theta)$-Bessel model}
   for $\pi$ is a subspace $\mathcal{B}_{\Lambda,\theta}(\pi)$ of
   $\mathcal{S}(\Lambda,\theta)$, invariant under right translation by
   $\mathfrak{g}$ and $K$, such that the $(\mathfrak{g},K)$-module thus defined
   is isomorphic to $(\pi,V)$.
  \item Let $(\pi,V)$ be a $(\mathfrak{g}^1,K^1)$-module.
   Let $\Lambda^1$ be a character of $T^1(\R)_0$. A \emph{$(\Lambda^1,\theta)$-Bessel model}
   for $\pi$ is a subspace $\mathcal{B}_{\Lambda,\theta}(\pi)$ of
   $\mathcal{S}^1(\Lambda^1,\theta)$, invariant under right translation by
   $\mathfrak{g}^1$ and $K^1$, such that the $(\mathfrak{g}^1,K^1)$-module thus defined
   is isomorphic to $(\pi,V)$.
 \end{enumerate}
\end{definition}
Our goal in the following is to prove uniqueness and give precise conditions
for existence of Bessel models for the lowest and highest
weight representations $\mathcal{E}_s(l,l')$ resp.\ $\mathcal{E}(l,l')$ described above.
One necessary condition is obvious: Since the Bessel subgroup $R(\R)$ contains the
center of $\GSp(4,\R)$ (resp.\ $\SSp(4,\R)$), the character $\Lambda$, restricted
to the center, has to coincide with the central character of the representation.
\subsubsection*{Change of models}
Let $S$ be a non-degenerate real symmetric matrix, as above.
Let $A\in\GL(2,\R)$ and $\alpha\in\R^\times$ be arbitrary, and define
$S'=\alpha\,^t\!ASA$. Then we have the two characters
$$
 \theta(\mat{1}{X}{}{1})=e^{2\pi i\,\tr(SX)},\qquad
 \theta'(\mat{1}{X}{}{1})=e^{2\pi i\,\tr(S'X)}.
$$
With $T(\R)=\{g\in\GL(2,\R):\:^tgSg=\det(g)S\}$ and
$T'(\R)=\{g\in\GL(2,\R):\:^tgS'g=\det(g)S'\}$, there is an isomorphism
\begin{align*}
 T'(\R)&\stackrel{\sim}{\longrightarrow}T(\R)\\
 g&\longmapsto AgA^{-1}.
\end{align*}
Let $\Lambda$ be a character of $T(\R)_0$, and let
$\Lambda'$ be the character of $T'(\R)_0$ corresponding to $\Lambda$, i.e.,
$\Lambda'(g)=\Lambda(AgA^{-1})$ for $g\in T'(\R)_0$.
Assume that $\mathcal{B}_{\Lambda,\theta}(\pi)$ is a Bessel model for
a $(\mathfrak{g},K)$-module $(\pi,V)$. Then it is easy to check
that, for $B\in\mathcal{B}(\Lambda,\theta)$, the function
$$
 B'(g)=B(\mat{A}{}{}{\alpha^{-1}\,^t\!A^{-1}}g),\qquad g\in\GSp(4,\R),
$$
satisfies $B'(tug)=\Lambda'(t)\theta'(u)B'(g)$ for $t\in T'(\R)_0$ and
$u\in U(\R)$. Hence, the map $B\mapsto B'$ provides a $(\mathfrak{g},K)$-isomorphism of
$\mathcal{B}_{\Lambda,\theta}(\pi)$ with a $(\Lambda',\theta')$-Bessel model
$\mathcal{B}_{\Lambda',\theta'}(\pi)$.
It follows that we need to prove existence and uniqueness of Bessel models
only for a class of representatives for quadratic forms $S$ under the operation
$S\mapsto\alpha\,^tASA$, $A\in\GL(2,\R)$, $\alpha\in\R^\times$.
There are only two such classes, represented by $S=\mat{1}{}{}{1}$ and $S=\mat{1}{}{}{-1}$. We call a Bessel model corresponding to $S=\mat{1}{}{}{1}$ a {\it non-split Bessel model} and a Bessel model corresponding to $S=\mat{1}{}{}{-1}$ a {\it split Bessel model}.
\nl
Similar considerations apply to Bessel models for $\SSp(4,\R)$. In this case
it is enough to prove existence and uniqueness of Bessel models
for a class of representatives for quadratic forms $S$ under the operation
$S\mapsto\,^tASA$, $A\in\GL(2,\R)$.
There are three such classes, represented by $S=\pm\mat{1}{}{}{1}$ and $S=\mat{1}{}{}{-1}$.
\nl
{\bf Remark:} Let $(\pi,V)$ be a $(\mathfrak{g}^1,K^1)$-module, and let $(\pi^\epsilon,V)$ be the
$(\mathfrak{g}^1,K^1)$-module defined in Sect.\ \ref{SpGSpsec}. Assume that $V$ is a
$(\Lambda^1,\theta)$-Bessel model for $\pi$. For $B\in V$, let
$\tilde B(g)=B(\epsilon g\epsilon^{-1})$, and let $\tilde V=\{\tilde B:\:B\in V\}$.
Then $\tilde V$ is a $(\mathfrak{g}^1,K^1)$-module under right translation, realizing
$\pi^\epsilon$. Evidently, $\tilde V$ is a $(\Lambda^1,\theta^{-1})$-Bessel model
for $\pi^\epsilon$. Hence,
\begin{equation}\label{pipiepsilonobviouseq}
 \text{$\pi$ has a $(\Lambda^1,\theta)$-Bessel model}\quad\Longleftrightarrow\quad
 \text{$\pi^\epsilon$ has a $(\Lambda^1,\theta^{-1})$-Bessel model}.
\end{equation}
\subsubsection*{Behavior under twisting}
Assume that $\mathcal{B}_{\Lambda,\theta}(\pi)$ is a $(\Lambda,\theta)$-Bessel model
for the representation $\pi$ of $\GSp(4,\R)$. Let $\chi$ be a character of $\R^\times$.
We attach to every $B\in\mathcal{B}_{\Lambda,\theta}(\pi)$
the function $\tilde B(g):=\chi(\mu_2(g))B(g)$, where $\mu_2$ is the multiplier
homomorphism. Let $\tilde V$ be the space of all functions $\tilde B$, where $B$
runs through $\mathcal{B}_{\Lambda,\theta}(\pi)$. Then right translation on $\tilde V$
defines a representation of $\GSp(4,\R)$ isomorphic to the twisted representation
$(\chi\otimes\pi)(g):=\chi(\mu_2(g))\pi(g)$. Each $\tilde B\in \tilde V$ satisfies
$$
 \tilde B(tug)=\chi(\det(t))\Lambda(t)\theta(u)\tilde B(g)\qquad
 \text{for }t\in T(\R)_0,\;u\in U(\R),\;g\in\GSp(4,\R).
$$
Here, $\det(t)$ is the determinant of $t$ considered as an element of $\GL(2,\R)$.
Since the multiplier function is slowly increasing, it follows that $\tilde V$
provides a $((\chi\circ\det)\Lambda,\theta)$-Bessel model for the twisted
representation. Hence,
\begin{equation}\label{twistedmodelseq}
 \text{$\pi$ has a $(\Lambda,\theta)$-Bessel model}\qquad\Longleftrightarrow\qquad
 \text{$\chi\otimes\pi$ has a $((\chi\circ\det)\Lambda,\theta)$-Bessel model}.
\end{equation}
Taking (\ref{TR0T1R0eq}) into account, it follows that in proving uniqueness and
existence of Bessel models, we may assume, whenever convenient, that the character
$\Lambda$ and the central character of $\pi$ are trivial on $\R_{>0}$.
\subsubsection*{Relating Bessel models for $\SSp(4,\R)$ and $\GSp(4,\R)$}
{\bf \underline{$\GSp(4)$ to $\SSp(4)$}:} Let $(\pi,V)$ be a given $(\mathfrak{g},K)$-module, and assume that
$V=\mathcal{B}_{\Lambda,\theta}(\pi)$ is a $(\Lambda,\theta)$-Bessel model.
Assume further that $\pi$ is irreducible, and that upon restriction to
$\SSp(4,\R)$ we have $V=V_1\oplus V_2$ with two non-isomorphic, irreducible
$(\mathfrak{g}^1,K^1)$-modules $(\pi_1,V_1)$ and $(\pi_2,V_2)$; see Sect.\ \ref{SpGSpsec}.
For $i=1,2$ let $\tilde V_i$ be the space of functions obtained by restricting
each function in $V_i$ to $\SSp(4,\R)$. The surjective map $V_i\rightarrow\tilde V_i$
given by restriction is obviously a $(\mathfrak{g}^1,K^1)$-map, and since $V_i$ is
irreducible, this map is either zero or an isomorphism. In case it is an isomorphism,
the space $\tilde V_i$ is a $(\Lambda^1,\theta)$-Bessel model for $\pi_i$.
It is clear that not all functions in $V$ can be supported on the non-identity
component of $\GSp(4,\R)$, so that at least one of the maps $V_i\rightarrow\tilde V_i$
must be non-zero. Hence, at least one of $\pi_1$ or $\pi_2$ admits a
$(\Lambda^1,\theta)$-Bessel model.
\nl
{\bf \underline{$\SSp(4)$ to $\GSp(4)$}:} Conversely, let $(\pi,V)$ be a given $(\mathfrak{g}^1,K^1)$-module for which
$\pi\not\cong\pi^\epsilon$, and assume that
$V=\mathcal{B}_{\Lambda^1,\theta}(\pi)$ is a $(\Lambda^1,\theta)$-Bessel model.
Let $\Lambda$ be any extension of $\Lambda^1$ to a character of $T(\R)_0$; see
(\ref{TR0T1R0eq}). Given $B\in V$, we extend $B$ to a function on $\SSp(4,\R)^\pm$
by setting it equal to zero on elements of negative multiplier, and then
further to a function on $\GSp(4,\R)=\R_{>0}\times\SSp(4,\R)^\pm$ which has
the $(\Lambda,\theta)$-Bessel transformation property;
see the discussion before Definition \ref{besselmodeldef}. Let $\tilde V$
be the space of functions thus obtained, and let $\tilde V^\epsilon$ be the
space of functions $\GSp(4,\R)\ni g\mapsto B(g\epsilon)$ for $B\in\tilde V$.
The spaces $\tilde V$ and $\tilde V^\epsilon$ have zero intersection, since
the functions in these spaces are supported on different connected components
of $\GSp(4,\R)$. The direct sum $\tilde V\oplus\tilde V^\epsilon$ is a
$(\mathfrak{g},K)$-module under right translation. It is a model for the
irreducible $(\mathfrak{g},K)$-module $I_s(\pi)$ considered in Sect.\ \ref{SpGSpsec},
where $s$ is the complex number defining the extension of $\Lambda^1$ to $\Lambda$.
Clearly, $\tilde V\oplus\tilde V^\epsilon$ is a $(\Lambda,\theta)$-Bessel model
for $I_s(\pi)$.
\nl
We just proved that a $(\Lambda^1,\theta)$-Bessel model for $\pi$ leads to a
$(\Lambda,\theta)$-Bessel model for $I_s(\pi)$, and clearly two different models
for $\pi$ would lead to two different models for $I_s(\pi)$.
Further below we will prove uniqueness of Bessel models for $\GSp(4,\R)$, which
therefore implies uniqueness of Bessel models for $\SSp(4,\R)$. It also shows that
$\pi$ and $\pi^\epsilon$ cannot both have a $(\Lambda^1,\theta)$-Bessel model, since
the above construction would lead to two different $(\Lambda,\theta)$-Bessel models
for $I_s(\pi)=I_s(\pi^\epsilon)$. We summarize:
\begin{proposition}\label{Sp4GSp4prop}
 Let $\Lambda$ be a character of $T(\R)_0$, let $\Lambda^1$ be its restriction
 to $T^1(\R)_0$, and let $\Lambda\big|_{\R_{>0}}$ be given by $a\mapsto a^s$ with
 $s\in\C$; see (\ref{TR0T1R0eq}). Let $(\pi,V)$ be a $(\mathfrak{g}^1,K^1)$-module
 for which $\pi\not\cong\pi^\epsilon$. Then the following are equivalent.
 \begin{enumerate}
  \item One of $\pi$ or $\pi^\epsilon$ has a $(\Lambda^1,\theta)$-Bessel model.
  \item Exactly one of $\pi$ or $\pi^\epsilon$ has a $(\Lambda^1,\theta)$-Bessel model.
  \item $I_s(\pi)$ has a $(\Lambda,\theta)$-Bessel model.
 \end{enumerate}
\end{proposition}
\section{Non-split Bessel models}\label{nonsplitsec}
In this section we investigate the existence and uniqueness of non-split Bessel models
for the lowest weight representations of $\GSp(4,\R)$ and
$\SSp(4,\R)$. We shall work with $\GSp(4,\R)$ and use the discussion preceding
Proposition \ref{Sp4GSp4prop} to obtain results for $\SSp(4,\R)$.
As explained in Sect.\ \ref{besselmodelssec}, we may throughout assume that
$$
 S=\mat{1}{}{}{1}.
$$
\subsection{Double coset decomposition}\label{nonsplitdoublecosetsec}
In this section we will derive representatives for the double coset space
$R(\R)\backslash\GSp(4,\R)/K^1$, where $R(\R)=T(\R)_0U(\R)$.
Recall from Sect.\ \ref{besselsubgroupssec} that
\begin{equation}\label{non-split-T(R)-gp}
 T(\R) = \{\mat{x}{y}{-y}{x}:\:x,y\in\R,\:x^2+y^2>0\}\cong\C^\times
\end{equation}
is a connected group.
The subgroup $T^1(\R)=T(\R)\cap\SL(2,\R)$ corresponds to the unit circle, and we have
\begin{equation}\label{TRTinfty1eq}
 T(\R)=\{\mat{\gamma}{}{}{\gamma}:\:\gamma>0\}\cdot T^1(\R).
\end{equation}
By the Cartan decomposition,
\begin{equation}\label{cartan1eq}
 \GL(2,\R)^+=\SO(2)\cdot\{\mat{\zeta_1}{}{}{\zeta_2}:\:\zeta_1\geq\zeta_2>0,\}\cdot\SO(2).
\end{equation}
Therefore,
\begin{align}\label{GL2RplusBesseldecompeq}
 \GL(2,\R)^+&=T^1(\R)\cdot\{\mat{\sqrt{\zeta_1\zeta_2}}{}{}{\sqrt{\zeta_1\zeta_2}}
 \mat{\sqrt{\zeta_1/\zeta_2}}{}{}{\sqrt{\zeta_2/\zeta_1}}
  :\:\zeta_1\geq\zeta_2>0\}\cdot\SO(2)\nonumber\\
 &=T(\R)\cdot\{\mat{\zeta}{}{}{\zeta^{-1}}:\:\zeta\geq1\}\cdot\SO(2).
\end{align}
Using this and the Iwasawa decomposition, it is not hard to see that
\begin{equation}\label{HRBesseldecompeq}
 \GSp(4,\R)=R(\R)\cdot\big\{\begin{bmatrix}\lambda \mat{\zeta}{}{}{\zeta^{-1}}&\\
  &\mat{\zeta^{-1}}{}{}{\zeta}\end{bmatrix}:\:
  \lambda\in\R^\times,\,\zeta\geq1\big\}\cdot K^1;
\end{equation}
see (4.7) of \cite{Fu}.
Here, $R(\R)=T(\R)U(\R)$ is the Bessel subgroup defined in (\ref{besselsubgroupdefeq}).
One can check that all the double cosets in (\ref{HRBesseldecompeq}) are disjoint.
Recall the coordinates (\ref{Kcoordeq}) in a neighborhood of the identity of $K^1$.
In the following we let
\begin{equation}\label{helementdefeq}
 h(\lambda,\zeta,\varphi_1,\varphi_2):=
 \begin{bmatrix}\lambda \mat{\zeta}{}{}{\zeta^{-1}}&\\
  &\mat{\zeta^{-1}}{}{}{\zeta}\end{bmatrix}r_1(\varphi_1)r_2(\varphi_2)
\end{equation}
for $\lambda,\zeta\in\R^\times$ and $\varphi_1,\varphi_2\in\R$.
\subsection{Differential operators}\label{diffopsecnonsplit}
In this section we will derive explicit formulas for the differential operators
given by elements of the complexified Lie algebra $\mathfrak{g}_\C$ on
the functions in a non-split Bessel model.
Assume that $\mathcal{B}_{\Lambda,\theta}(\pi)$ is a Bessel model for the
$(\mathfrak{g},K)$-module $(\pi,V)$.
For any $B\in\mathcal{B}_{\Lambda,\theta}(\pi)$ we define a function
$f=f_B$ on $\R^\times\times\R_{>0}\times\R\times\R$ by
\begin{equation}\label{fdefeq}
 f(\lambda,\zeta,\varphi_1,\varphi_2)=B(h(\lambda,\zeta,\varphi_1,\varphi_2)).
\end{equation}
It follows from (\ref{HRBesseldecompeq}) and the $K$-finiteness of $B$
that if $B$ has weight $(l,l')$, then
$B$ is determined by $f$. If $L$ denotes one of the
operators $N_\pm,X_\pm,P_{0\pm},P_{1\pm}$, then $L.B$ will
be determined by the associated function $f_{L.B}$.
\nl
We first have to compute the action of the non-complexified
Lie algebra $\mathfrak g$. If $L\in\mathfrak g$, then by definition
$$
 (L.B)(h(\lambda,\zeta,\varphi_1,\varphi_2))=
 \frac d{dt}\Big|_0B\big(h(\lambda,\zeta,\varphi_1,\varphi_2)\exp(tL)\big).
$$
Now, at least for small values of $t$, we can decompose the argument according to
(\ref{HRBesseldecompeq}),
\begin{align}\label{exptLiwasawaeq}
 h(\lambda,\zeta,\varphi_1,\varphi_2)\exp(tL)&=\mat{g(t)}{}{}{\det(g(t))\,^tg(t)^{-1}}
 \begin{bmatrix}1&&x(t)&y(t)\\&1&y(t)&z(t)\\&&1\\&&&1\end{bmatrix}\nonumber\\
 &\hspace{20ex}h(\lambda(t),\zeta(t),\varphi_1(t),\varphi_2(t))
 \,r_3(\varphi_3(t))r_4(\varphi_4(t)).
\end{align}
Here, $g(t)\in T(\R)$, and $x(t)$ etc.\ are smooth functions in a neighborhood of $0$
satisfying
\begin{align*}
 &x(0)=y(0)=z(0)=\varphi_3(0)=\varphi_4(0)=0,\\
 &\lambda(0)=\lambda,\quad\zeta(0)=\zeta,
 \quad\varphi_1(0)=\varphi_1,\quad\varphi_2(0)=\varphi_2.
\end{align*}
According to (\ref{TRTinfty1eq}), we can write
\begin{equation}\label{gtdecompeq}
 g(t)=\gamma(t)\mat{\cos(\delta(t))}{\sin(\delta(t))}{-\sin(\delta(t))}{\cos(\delta(t))}
\end{equation}
with smooth functions $\gamma(t)$ and $\delta(t)$ such that $\gamma(0)=1$ and
$\delta(0)=0$. The character $\Lambda$ of $T(\R)$ is of the form
\begin{equation}\label{Lambdaexpliciteq}
 \Lambda(\gamma \mat{\cos(\delta)}{\sin(\delta)}{-\sin(\delta)}{\cos(\delta)})
 =\gamma^se^{im\delta},\qquad\gamma>0,\;\delta\in\R,
\end{equation}
with some $s\in\C$ and $m\in\Z$.
It follows that
\begin{align}\label{Loperationgeneraleq}
 &(L.B)(h(\lambda,\zeta,\varphi_1,\varphi_2))
  =\frac d{dt}\Big|_0\Big(\Lambda(g(t))\theta(\begin{bmatrix}1&&x(t)&y(t)\\&1&y(t)&z(t)\\
  &&1\\&&&1\end{bmatrix})e^{i(l\varphi_3(t)+l'\varphi_4(t))}
  B\big(h(\lambda(t),\zeta(t),\varphi_1(t),\varphi_2(t))\big)\Big)\nonumber\\
 &\;=\frac d{dt}\Big|_0\Big(\gamma(t)^se^{im\delta(t)}e^{2\pi i(ax(t)+by(t)+cz(t))}
  e^{i(l\varphi_3(t)+l'\varphi_4(t))}
  f\big(\lambda(t),\zeta(t),\varphi_1(t),\varphi_2(t)\big)\Big)\nonumber\\
 &\;=\Big(s\gamma'(0)+i(m\delta'(0)+l\varphi_3'(0)+l'\varphi_4'(0))
  +2\pi i(ax'(0)+by'(0)+cz'(0))\Big)f(\lambda,\zeta,\varphi_1,\varphi_2)\nonumber\\
 &\quad+
  \lambda'(0)\frac{\partial f}{\partial\lambda}(\lambda,\zeta,\varphi_1,\varphi_2)
  +\zeta'(0)\frac{\partial f}{\partial\zeta}(\lambda,\zeta,\varphi_1,\varphi_2)
  +\varphi_1'(0)\frac{\partial f}{\partial\varphi_1}(\lambda,\zeta,\varphi_1,\varphi_2)
  +\varphi_2'(0)\frac{\partial f}{\partial\varphi_2}(\lambda,\zeta,\varphi_1,\varphi_2).
\end{align}
Thus, what we need are the derivatives at $0$ of the auxiliary functions
$\gamma,\delta,\ldots$. To get these, we differentiate the
matrix equation (\ref{exptLiwasawaeq}) and put $t=0$. This yields
sixteen linear equations from which the desired derivatives can be
determined. The results are as follows.
\begin{enumerate}
\item Let $L=H_1=\begin{bmatrix}1&0&0&0\\0&0&0&0\\0&0&-1&0\\0&0&0&0\end{bmatrix}$. Then we have
    \begin{eqnarray}
    &&\gamma'(0)=-\frac 12 \cos(2 \varphi_2), \,\, \delta'(0) = \frac{\zeta^2 \sin(2 \varphi_1)}{\zeta^4-1},\,\, \lambda'(0)=\lambda\cos(2\varphi_2), \,\, \zeta'(0)=\frac 12 \zeta\cos(2\varphi_1), \nonumber\\
    && x'(0)=-\zeta^2 \lambda\sin(2\varphi_1)\sin(2\varphi_2), \,y'(0)=-\lambda \cos(2\varphi_1)\sin(2\varphi_2), \,z'(0)=\frac{\lambda \sin(2\varphi_1)\sin(2\varphi_2)}{\zeta^2}, \nonumber\\
    &&\varphi_1'(0)=\frac{(1+\zeta^4)\sin(2\varphi_1)}{2(1-\zeta^4)}, \,\, \varphi_2'(0)=\frac 12 \sin(2\varphi_2), \,\, \varphi_3'(0)=0, \,\, \varphi_4'(0)=0 \label{H1constants}
    \end{eqnarray}

\item Let $L=H_2=\begin{bmatrix}0&0&0&0\\0&1&0&0\\0&0&0&0
  \\0&0&0&-1\end{bmatrix}$. Then we have
  \begin{eqnarray}
    &&\gamma'(0)=-\frac 12 \cos(2 \varphi_2), \,\, \delta'(0) = -\frac{\zeta^2 \sin(2 \varphi_1)}{\zeta^4-1},\,\, \lambda'(0)=\lambda\cos(2\varphi_2), \,\, \zeta'(0)=-\frac 12 \zeta\cos(2\varphi_1), \nonumber\\
    && x'(0)=-\zeta^2 \lambda\sin(2\varphi_1)\sin(2\varphi_2), \,y'(0)=-\lambda \cos(2\varphi_1)\sin(2\varphi_2), \,z'(0)=\frac{\lambda \sin(2\varphi_1)\sin(2\varphi_2)}{\zeta^2}, \nonumber\\
    &&\varphi_1'(0)=\frac{(1+\zeta^4)\sin(2\varphi_1)}{2(\zeta^4-1)}, \,\, \varphi_2'(0)=\frac 12 \sin(2\varphi_2), \,\, \varphi_3'(0)=0, \,\, \varphi_4'(0)=0 \label{H2constants}
    \end{eqnarray}

\item Let $L=F=\begin{bmatrix}0&0&1&0\\0&0&0&0\\0&0&0&0\\0&0&0&0\end{bmatrix}$. Then we have
\begin{eqnarray}
    &&\gamma'(0)=0, \,\, \delta'(0) = \frac{\zeta^2\cos(2\varphi_1)\sin(2\varphi_2)}{2(1-\zeta^4)},\,\, \lambda'(0)=0, \,\, \zeta'(0)=\frac{\zeta}4\sin(2\varphi_1)\sin(2\varphi_2), \nonumber\\
    && x'(0)=\frac 12\zeta^2 \lambda(\cos(2\varphi_1)+\cos(2\varphi_2)), \,y'(0)=-\frac{\lambda}2 \sin(2\varphi_1), \,z'(0)=\frac{\lambda (\cos(2\varphi_2)-\cos(2\varphi_1))}{2\zeta^2}, \nonumber\\
    &&\varphi_1'(0)=\frac 14 \tan(2\varphi_2) +\frac{\zeta^4+1}{4(\zeta^4-1)}\cos(2\varphi_1)\sin(2\varphi_2), \,\, \varphi_2'(0)=0, \,\, \varphi_3'(0)=-\frac{\sin(\varphi_2)^4}{\cos(2\varphi_2)}, \nonumber\\
     &&\varphi_4'(0)=\frac 14 \sin(2\varphi_2)\tan(2\varphi_2)\label{Fconstants}
    \end{eqnarray}

\item Let $L=G=\begin{bmatrix}0&0&0&0\\0&0&0&0\\1&0&0&0
  \\0&0&0&0\end{bmatrix}$. Then we have
  \begin{eqnarray}
    &&\gamma'(0)=0, \,\, \delta'(0) = \frac{\zeta^2\cos(2\varphi_1)\sin(2\varphi_2)}{2(1-\zeta^4)},\,\, \lambda'(0)=0, \,\, \zeta'(0)=\frac{\zeta}4\sin(2\varphi_1)\sin(2\varphi_2), \nonumber\\
    && x'(0)=\frac 12\zeta^2 \lambda(\cos(2\varphi_1)+\cos(2\varphi_2)), \,y'(0)=-\frac{\lambda}2 \sin(2\varphi_1), \,z'(0)=\frac{\lambda (\cos(2\varphi_2)-\cos(2\varphi_1))}{2\zeta^2}, \nonumber\\
    &&\varphi_1'(0)=\frac 14 \tan(2\varphi_2) +\frac{\zeta^4+1}{4(\zeta^4-1)}\cos(2\varphi_1)\sin(2\varphi_2), \,\, \varphi_2'(0)=0, \,\, \varphi_3'(0)=-\frac{\cos(\varphi_2)^4}{\cos(2\varphi_2)}, \nonumber\\
     &&\varphi_4'(0)=\frac 14 \sin(2\varphi_2)\tan(2\varphi_2)\label{Gconstants}
    \end{eqnarray}

\item Let $L=R=\begin{bmatrix}0&0&0&0\\0&0&0&1\\0&0&0&0\\0&0&0&0\end{bmatrix}$. Then we have
    \begin{eqnarray}
    &&\gamma'(0)=0, \,\, \delta'(0) = \frac{\zeta^2\cos(2\varphi_1)\sin(2\varphi_2)}{2(1-\zeta^4)},\,\, \lambda'(0)=0, \,\, \zeta'(0)=\frac{\zeta}4\sin(2\varphi_1)\sin(2\varphi_2), \nonumber\\
    && x'(0)=\frac 12\zeta^2 \lambda(\cos(2\varphi_2)-\cos(2\varphi_1)), \,y'(0)=\frac{\lambda}2 \sin(2\varphi_1), \,z'(0)=\frac{\lambda (\cos(2\varphi_2)+\cos(2\varphi_1))}{2\zeta^2}, \nonumber\\
    &&\varphi_1'(0)=-\frac 14 \tan(2\varphi_2) +\frac{\zeta^4+1}{4(\zeta^4-1)}\cos(2\varphi_1)\sin(2\varphi_2), \,\, \varphi_2'(0)=0, \,\, \varphi_3'(0)=\frac{\sin(2\varphi_2)^2}{4\cos(2\varphi_2)}, \nonumber\\
     &&\varphi_4'(0)=-\frac{\sin(\varphi_2)^4}{\cos(2\varphi_2)}\label{Rconstants}
    \end{eqnarray}

\item Let $L=R'=\begin{bmatrix}0&0&0&0\\0&0&0&0\\0&0&0&0
  \\0&1&0&0\end{bmatrix}$. Then we have
  \begin{eqnarray}
    &&\gamma'(0)=0, \,\, \delta'(0) = \frac{\zeta^2\cos(2\varphi_1)\sin(2\varphi_2)}{2(1-\zeta^4)},\,\, \lambda'(0)=0, \,\, \zeta'(0)=\frac{\zeta}4\sin(2\varphi_1)\sin(2\varphi_2), \nonumber\\
    && x'(0)=\frac 12\zeta^2 \lambda(\cos(2\varphi_2)-\cos(2\varphi_1)), \,y'(0)=\frac{\lambda}2 \sin(2\varphi_1), \,z'(0)=\frac{\lambda (\cos(2\varphi_2)+\cos(2\varphi_1))}{2\zeta^2}, \nonumber\\
    &&\varphi_1'(0)=-\frac 14 \tan(2\varphi_2) +\frac{\zeta^4+1}{4(\zeta^4-1)}\cos(2\varphi_1)\sin(2\varphi_2), \,\, \varphi_2'(0)=0, \,\, \varphi_3'(0)=\frac{\sin(2\varphi_2)^2}{4\cos(2\varphi_2)}, \nonumber\\
     &&\varphi_4'(0)=-\frac{\cos(\varphi_2)^4}{\cos(2\varphi_2)}\label{R'constants}
    \end{eqnarray}

\item Let $L=P=\begin{bmatrix}0&0&0&0\\1&0&0&0\\0&0&0&-1\\0&0&0&0\end{bmatrix}$. Then we have
    \begin{eqnarray}
    &&\gamma'(0)=0, \,\, \delta'(0)=\frac{\zeta^2\cos(2\varphi_1)\cos(2\varphi_2)}{1-\zeta^4}, \,\, \lambda'(0)=0, \,\, \zeta'(0)=\frac 12\zeta \cos(2\varphi_2)\sin(2\varphi_1) \nonumber\\
    &&x'(0)=-\zeta^2 \lambda\sin(2\varphi_2), \,\, y'(0)=0, \,\, z'(0)=-\frac{\lambda \sin(2\varphi_2)}{\zeta^2} \nonumber \\
    &&\varphi_1'(0)=-\frac{1}{2\cos(2\varphi_2)}+\frac{\zeta^4+1}{2(\zeta^4-1)}\cos(2\varphi_1)\cos(2\varphi_2), \,\, \varphi_2'(0)=0, \,\, \varphi_3'(0)=\tan(2\varphi_2)\cos(\varphi_2)^2, \nonumber\\
    &&\varphi_4'(0)=-\tan(2\varphi_2)\sin(\varphi_2)^2 \label{Pconstants}
    \end{eqnarray}

\item Let $L=P'=\begin{bmatrix}0&1&0&0\\0&0&0&0\\0&0&0&0
  \\0&0&-1&0\end{bmatrix}$. Then we have
  \begin{eqnarray}
    &&\gamma'(0)=0, \,\, \delta'(0)=\frac{\zeta^2\cos(2\varphi_1)\cos(2\varphi_2)}{1-\zeta^4}, \,\, \lambda'(0)=0, \,\, \zeta'(0)=\frac 12\zeta \cos(2\varphi_2)\sin(2\varphi_1) \nonumber\\
    &&x'(0)=-\zeta^2 \lambda\sin(2\varphi_2), \,\, y'(0)=0, \,\, z'(0)=-\frac{\lambda \sin(2\varphi_2)}{\zeta^2} \nonumber \\
    &&\varphi_1'(0)=\frac{1}{2\cos(2\varphi_2)}+\frac{\zeta^4+1}{2(\zeta^4-1)}\cos(2\varphi_1)\cos(2\varphi_2), \,\, \varphi_2'(0)=0, \,\, \varphi_3'(0)=-\tan(2\varphi_2)\sin(\varphi_2)^2, \nonumber\\
    &&\varphi_4'(0)=\tan(2\varphi_2)\cos(\varphi_2)^2 \label{P'constants}
    \end{eqnarray}

\item Let $L=Q=\begin{bmatrix}0&0&0&1\\0&0&1&0\\0&0&0&0\\0&0&0&0\end{bmatrix}$. Then we have
    \begin{eqnarray}
    &&\gamma'(0)=-\frac{\sin(2\varphi_2)}2, \,\, \delta'(0)=0, \,\, \lambda'(0)= \lambda \sin(2\varphi_2), \,\, \zeta'(0)=0 \nonumber\\
    &&x'(0)=\zeta^2\lambda\cos(2\varphi_2)\sin(2\varphi_1), \,\, y'(0)=\lambda \cos(2\varphi_2)\cos(2\varphi_1), \,\, z'(0)=-\frac{\lambda\cos(2\varphi_2)\sin(2\varphi_1)}{\zeta^2} \nonumber\\
    &&\varphi_1'(0)=0, \,\, \varphi_2'(0)=\sin(\varphi_2)^2, \,\, \varphi_3'(0)=0, \,\, \varphi_4'(0)=0 \label{Qconstant}
    \end{eqnarray}

\item Let $L=Q'=\begin{bmatrix}0&0&0&0\\0&0&0&0\\0&1&0&0
  \\1&0&0&0\end{bmatrix}$. Then we have
  \begin{eqnarray}
    &&\gamma'(0)=-\frac{\sin(2\varphi_2)}2, \,\, \delta'(0)=0, \,\, \lambda'(0)= \lambda \sin(2\varphi_2), \,\, \zeta'(0)=0 \nonumber\\
    &&x'(0)=\zeta^2\lambda\cos(2\varphi_2)\sin(2\varphi_1), \,\, y'(0)=\lambda \cos(2\varphi_2)\cos(2\varphi_1), \,\, z'(0)=-\frac{\lambda\cos(2\varphi_2)\sin(2\varphi_1)}{\zeta^2} \nonumber\\
    &&\varphi_1'(0)=0, \,\, \varphi_2'(0)=-\cos(\varphi_2)^2, \,\, \varphi_3'(0)=0, \,\, \varphi_4'(0)=0 \label{Q'constant}
    \end{eqnarray}
\end{enumerate}
Using these coefficients, we obtain from (\ref{Loperationgeneraleq}) the following
formulas for the action of the elements of the complexified Lie algebra.
\begin{align}
 Z.B&=lB,\label{Z-formula}\\
 Z'.B&=l'B, \label{Z'-formula}\\
 (N_\pm.B)(h(\lambda,\zeta,\varphi_1,\varphi_2))
 &=\frac i2\tan(2\varphi_2)(l'-l)f(h(\lambda,\zeta,\varphi_1,\varphi_2))
  +\frac1{2\cos(2\varphi_2)}\frac{\partial f}{\partial\varphi_1}
   (h(\lambda,\zeta,\varphi_1,\varphi_2)) \nonumber \\
  &\quad\mp\frac i2\frac{\partial f}{\partial\varphi_2}
   (h(\lambda,\zeta,\varphi_1,\varphi_2)), \label{N-formula}\\
 (X_\pm.B)(h(\lambda,\zeta,\varphi_1,\varphi_2))
  &\;=\Big(-\frac s4\cos(2\varphi_2)+\frac m2\frac1{\zeta^2-\zeta^{-2}}
  \big(i\sin(2\varphi_1)\pm\cos(2\varphi_1)\sin(2\varphi_2)\big) \nonumber\\
  &\quad\pm\frac{l}2\frac{\sin(\varphi_2)^4+\cos(\varphi_2)^4}{\cos(2\varphi_2)}
   \mp\frac{l'}4\frac{\sin(2\varphi_2)^2}{\cos(2\varphi_2)}
   -\pi i\lambda(\zeta^2-\zeta^{-2})\sin(2\varphi_1)\sin(2\varphi_2) \nonumber\\
  &\quad\mp\pi\lambda
   (\zeta^2-\zeta^{-2})\cos(2\varphi_1)\mp\pi\lambda(\zeta^2+\zeta^{-2})\cos(2\varphi_2)
   \Big)f(\lambda,\zeta,\varphi_1,\varphi_2)\nonumber\\
 &\quad+\frac12\cos(2\varphi_2)\lambda
  \frac{\partial f}{\partial\lambda}(\lambda,\zeta,\varphi_1,\varphi_2) \nonumber\\
 &\quad+\frac14\big(\cos(2\varphi_1)\pm i\sin(2\varphi_1)\sin(2\varphi_2)\big)\zeta
  \frac{\partial f}{\partial\zeta}(\lambda,\zeta,\varphi_1,\varphi_2) \nonumber\\
 &\quad+\Big(\frac{\zeta^2+\zeta^{-2}}{4(\zeta^2-\zeta^{-2})}
  \big(\!-\!\sin(2\varphi_1)\pm i\cos(2\varphi_1)
  \sin(2\varphi_2)\big)\pm\frac i4\tan(2\varphi_2)\Big)
  \frac{\partial f}{\partial\varphi_1}(\lambda,\zeta,\varphi_1,\varphi_2) \nonumber\\
 &\quad+\frac14\sin(2\varphi_2)\frac{\partial f}{\partial\varphi_2}
  (\lambda,\zeta,\varphi_1,\varphi_2),\label{X-formula}\\
 (P_{1\pm}.B)(h(\lambda,\zeta,\varphi_1,\varphi_2))
  &=\Big(\mp\frac{is}2\sin(2\varphi_2)-\frac{im}{\zeta^2-\zeta^{-2}}
   \cos(2\varphi_1)\cos(2\varphi_2)+\frac{i(l+l')}2\sin(2\varphi_2)\nonumber\\
  &\quad-2\pi i\lambda(\zeta^2+\zeta^{-2})\sin(2\varphi_2)
  \mp2\pi\lambda(\zeta^2-\zeta^{-2})\cos(2\varphi_2)\sin(2\varphi_1)
  \Big)f(\lambda,\zeta,\varphi_1,\varphi_2)\nonumber\\
 &\quad\pm i\sin(2\varphi_2)\lambda
  \frac{\partial f}{\partial\lambda}(\lambda,\zeta,\varphi_1,\varphi_2) \nonumber\\
 &\quad+\frac12\cos(2\varphi_2)\sin(2\varphi_1)\zeta
  \frac{\partial f}{\partial\zeta}(\lambda,\zeta,\varphi_1,\varphi_2) \nonumber\\
 &\quad+\frac{\zeta^2+\zeta^{-2}}{2(\zeta^2-\zeta^{-2})}\cos(2\varphi_1)\cos(2\varphi_2)
  \frac{\partial f}{\partial\varphi_1}(\lambda,\zeta,\varphi_1,\varphi_2) \nonumber\\
 &\quad\mp\frac i2\cos(2\varphi_2)\frac{\partial f}{\partial\varphi_2}
  (\lambda,\zeta,\varphi_1,\varphi_2), \label{P1-formula}\\
 (P_{0\pm}.B)(h(\lambda,\zeta,\varphi_1,\varphi_2))
  &=\Big(-\frac s4\cos(2\varphi_2)-\frac m2\frac1{\zeta^2-\zeta^{-2}}
   (i\sin(2\varphi_1)\mp \cos(2\varphi_1)\sin(2\varphi_2)) \nonumber\\
   &\quad\mp\frac l4\frac{\sin(2\varphi_2)^2}{\cos(2\varphi_2)}
   \pm\frac{l'}2\frac{\sin(\varphi_2)^4+\cos(\varphi_2)^4}{\cos(2\varphi_2)}
   -\pi i\lambda(\zeta^2-\zeta^{-2})\sin(2\varphi_1)\sin(2\varphi_2) \nonumber\\
  &\quad\mp\pi\lambda(\zeta^2+\zeta^{-2})\cos(2\varphi_2)
   \mp\pi\lambda(\zeta^{-2}-\zeta^2)\cos(2\varphi_1)
   \Big)f(\lambda,\zeta,\varphi_1,\varphi_2)\nonumber\\
 &\quad+\frac12\cos(2\varphi_2)\lambda
  \frac{\partial f}{\partial\lambda}(\lambda,\zeta,\varphi_1,\varphi_2) \nonumber\\
 &\quad+\frac14(-\cos(2\varphi_1)\pm i\sin(2\varphi_1)\sin(2\varphi_2))\zeta
  \frac{\partial f}{\partial\zeta}(\lambda,\zeta,\varphi_1,\varphi_2) \nonumber\\
 &\quad+\Big(\frac{\zeta^2+\zeta^{-2}}{4(\zeta^2-\zeta^{-2})}
  (\sin(2\varphi_1)\pm i\cos(2\varphi_1)\sin(2\varphi_2))\mp\frac i4\tan(2\varphi_2)\Big)
  \frac{\partial f}{\partial\varphi_1}(\lambda,\zeta,\varphi_1,\varphi_2) \nonumber\\
 &\quad+\frac 14\sin(2\varphi_2)
  \frac{\partial f}{\partial\varphi_2}(\lambda,\zeta,\varphi_1,\varphi_2).\label{P0-formula}
\end{align}
\subsection{Existence and uniqueness of Bessel functions}\label{hdsrnonsplitsec}
Let $\Lambda$ be the character of $\C^{\times}$ defined in (\ref{Lambdaexpliciteq}), i.e.,
\begin{equation}\label{Lambdaexpliciteq2}
 \Lambda(\gamma \mat{\cos(\delta)}{\sin(\delta)}{-\sin(\delta)}{\cos(\delta)})
 =\gamma^se^{im\delta},\qquad\gamma>0,\;\delta\in\R,
\end{equation}
with $s\in\C$ and $m\in\Z$.
Let $\pi=\mathcal{E}(l,l')$ be the lowest weight representation of
$\GSp(4,\R)$ with minimal $K$-type $(l,l')$.
Assume that $B$ is a highest weight vector in the minimal $K$-type in a
Bessel model for $\pi$ of type $(\Lambda,\theta)$. Hence, $B$ satisfies
the following conditions.
\begin{enumerate}
 \item[($\mathcal{S}1$)] $B$ is smooth and $K$-finite.
 \item[($\mathcal{S}2$)] $B(tug)=\Lambda(t)\theta(u)B(g)$ for all
  $t\in T(\R)$, $u\in U(\R)$, $g\in\GSp(4,\R)$.
 \item[($\mathcal{S}3$)] $B$ is slowly increasing.
 \item[($\mathcal{S}4$)] $Z.B=lB$ and $Z'.B = l'B$. Equivalently,
  $B(gr_3(\varphi_3)r_4(\varphi_4))=e^{i(l\varphi_3+l'\varphi_4)}B(g)$
  for all $\varphi_3,\varphi_4\in\R$, $g\in\GSp(4,\R)$.
 \item[($\mathcal{S}5$)] $N_+.B=X_-.B=P_{1-}.B=P_{0-}.B=0$.
\end{enumerate}
Conditions ($\mathcal{S}1$) -- ($\mathcal{S}3$) say that $B$ is an element
of the space $\mathcal{S}(\Lambda,\theta)$ defined in Sect.\ \ref{besselmodelssec}.
Condition ($\mathcal{S}4$) says that $B$ has weight $(l,l')$.
And condition ($\mathcal{S}5$) means that $B$ is a lowest weight vector.
\nl
Let $\mathcal{S}(\Lambda,\theta,l,l')$ be the space of all functions
$B:\:\GSp(4,\R)\rightarrow\C$ satisfying conditions ($\mathcal{S}1$) -- ($\mathcal{S}5$).
It follows from ($\mathcal{S}4$) that $B(-g)=(-1)^{l+l'}B(g)$. On the other hand, by ($\mathcal{S}2$) and
(\ref{Lambdaexpliciteq2}), $B(-g)=(-1)^mB(g)$. Hence, a necessary condition
for $\mathcal{S}(\Lambda,\theta,l,l')$ to be non-zero
is that $l+l'+m$ is an even integer.
\subsubsection*{Uniqueness in the neighborhood of identity}
As in (\ref{fdefeq}) we let
$f(\lambda,\zeta,\varphi_1,\varphi_2)=B(h(\lambda,\zeta,\varphi_1,\varphi_2))$.
From formulas (\ref{N-formula}) -- (\ref{P0-formula}) we find that,
in a neighborhood of the identity, the conditions in ($\mathcal{S}5$)
are equivalent to the following system of linear first-order PDEs.
\begin{align}
 \frac{\partial f}{\partial\lambda}&=\Big(\frac{l+l'+s}{2\lambda}-2\pi(\zeta^2+\zeta^{-2})
  \Big)f,\label{general1eq1}\\
 \frac{\partial f}{\partial\zeta}&=\Big(-4\pi\lambda(\zeta-\zeta^{-3})-
  \frac{2m\zeta(\cos(2\varphi_1)\sin(2\varphi_2)+i\sin(2\varphi_1))}
  {(\zeta^4-1)(\cos(2\varphi_1)+i\sin(2\varphi_1)\sin(2\varphi_2))
  +(\zeta^4+1)\cos(2\varphi_2)} \nonumber \\
 &\qquad+\frac{l-l'}{\zeta}\frac{(\zeta^4+1)(\cos(2\varphi_1)
  +i\sin(2\varphi_1)\sin(2\varphi_2))
  +(\zeta^4-1)\cos(2\varphi_2)}
  {(\zeta^4-1)(\cos(2\varphi_1)+i\sin(2\varphi_1)\sin(2\varphi_2))
  +(\zeta^4+1)\cos(2\varphi_2)}\Big)f,\label{general2eq1}\\
 \frac{\partial f}{\partial\varphi_1}&=\Big(
  \frac{2im\zeta^2\cos(2\varphi_2)-(l-l')(\zeta^4-1)(\sin(2\varphi_1)-
  i\cos(2\varphi_1)\sin(2\varphi_2))}
  {(\zeta^4-1)(\cos(2\varphi_1)+i\sin(2\varphi_1)\sin(2\varphi_2))
  +(\zeta^4+1)\cos(2\varphi_2)}\Big)f,\label{general3eq1}\\
 \frac{\partial f}{\partial\varphi_2}&=\Big(
  \frac{2m\zeta^2-(l-l')\big(
  (\zeta^4+1)\sin(2\varphi_2)-i(\zeta^4-1)\sin(2\varphi_1)\cos(2\varphi_2)\big)}
  {(\zeta^4-1)(\cos(2\varphi_1)+i\sin(2\varphi_1)\sin(2\varphi_2))
  +(\zeta^4+1)\cos(2\varphi_2)}\Big)f.\label{general4eq1}
\end{align}
From (\ref{general1eq1}) we conclude
\begin{equation}\label{general5eq1}\renewcommand{\arraystretch}{1.4}
 f(\lambda,\zeta,\varphi_1,\varphi_2)=\left\{\begin{array}{l@{\qquad\text{if }}l}
 c_1(\zeta,\varphi_1,\varphi_2)\lambda^{\frac{l+l'+s}2}
 e^{-2\pi\lambda(\zeta^2+\zeta^{-2})}&\lambda>0,\\
 c_2(\zeta,\varphi_1,\varphi_2)(-\lambda)^{\frac{l+l'+s}2}
 e^{-2\pi\lambda(\zeta^2+\zeta^{-2})}&\lambda<0,\end{array}\right.
\end{equation}
for some functions $c_1(\zeta,\varphi_1,\varphi_2)$ and
$c_2(\zeta,\varphi_1,\varphi_2)$. Since we are only interested in slowly
increasing solutions, we must have $c_2(\zeta,\varphi_1,\varphi_2)=0$;
see (\ref{slowlyincreasingeq2}). Substituting the solution for $\lambda>0$
into (\ref{general2eq1}), (\ref{general3eq1}) and (\ref{general4eq1}),
we obtain the system
\begin{align}
 \frac{\partial c_1}{\partial\zeta}&=\Big(
  -\frac{2m\zeta(\cos(2\varphi_1)\sin(2\varphi_2)+i\sin(2\varphi_1))}
  {(\zeta^4-1)(\cos(2\varphi_1)+i\sin(2\varphi_1)\sin(2\varphi_2))
  +(\zeta^4+1)\cos(2\varphi_2)} \nonumber \\
 &\qquad+\frac{l-l'}{\zeta}\frac{(\zeta^4+1)(\cos(2\varphi_1)
  +i\sin(2\varphi_1)\sin(2\varphi_2))+(\zeta^4-1)\cos(2\varphi_2)}
  {(\zeta^4-1)(\cos(2\varphi_1)+i\sin(2\varphi_1)\sin(2\varphi_2))
  +(\zeta^4+1)\cos(2\varphi_2)}\Big)c_1,\label{general6eq1}\\
 \frac{\partial c_1}{\partial\varphi_1}&=\Big(
  \frac{2im\zeta^2\cos(2\varphi_2)-(l-l')(\zeta^4-1)
  (\sin(2\varphi_1)-i\cos(2\varphi_1)\sin(2\varphi_2))}
  {(\zeta^4-1)(\cos(2\varphi_1)+i\sin(2\varphi_1)\sin(2\varphi_2))
  +(\zeta^4+1)\cos(2\varphi_2)}\Big)c_1,\label{general7eq1}\\
 \frac{\partial c_1}{\partial\varphi_2}&=\Big(
  \frac{2m\zeta^2-(l-l')\big((\zeta^4+1)\sin(2\varphi_2)-i(\zeta^4-1)
  \sin(2\varphi_1)\cos(2\varphi_2)\big)}
  {(\zeta^4-1)(\cos(2\varphi_1)+i\sin(2\varphi_1)\sin(2\varphi_2))
  +(\zeta^4+1)\cos(2\varphi_2)}\Big)c_1\label{general8eq1}
\end{align}
for the function $c_1$.
\begin{lemma}\label{holdiscserminKBessellemma}
 Let $l,l',m$ be integers such that $l+l'+m$ is even. Then the system
 (\ref{general6eq1}), (\ref{general7eq1}), (\ref{general8eq1}) has, up to scalars,
 the unique solution
 \begin{align}\label{holdiscserminKBessellemmaeq1}
  \nonumber c_1(\zeta,\varphi_1,\varphi_2)&=
  \Big(\zeta(\cos(\varphi_1)\cos(\varphi_2)+i\sin(\varphi_1)\sin(\varphi_2))
  +\zeta^{-1}(\cos(\varphi_1)\sin(\varphi_2)+i\sin(\varphi_1)\cos(\varphi_2))\Big)^{m}\\
  &\qquad\Big((\zeta^2-\zeta^{-2})(\cos(2\varphi_1)
  +i\sin(2\varphi_1)\sin(2\varphi_2))+(\zeta^2+\zeta^{-2})\cos(2\varphi_2)\Big)
  ^{\frac{l-l'-m}2}.
 \end{align}
 Alternatively,
 \begin{align}\label{holdiscserminKBessellemmaeq2}
  \nonumber c_1(\zeta,\varphi_1,\varphi_2)&=2^{-m}
  \Big(\zeta(\cos(\varphi_1)\cos(\varphi_2)+i\sin(\varphi_1)\sin(\varphi_2))
  -\zeta^{-1}(\cos(\varphi_1)\sin(\varphi_2)+i\sin(\varphi_1)\cos(\varphi_2))\Big)^{-m}\\
  &\qquad\Big((\zeta^2-\zeta^{-2})(\cos(2\varphi_1)
  +i\sin(2\varphi_1)\sin(2\varphi_2))+(\zeta^2+\zeta^{-2})\cos(2\varphi_2)\Big)
  ^{\frac{l-l'+m}2}.
 \end{align}
\end{lemma}
{\bf Proof:} A direct calculation shows that the two expressions are equal.
It is easily verified that (\ref{holdiscserminKBessellemmaeq1}) satisfies
(\ref{general6eq1}), (\ref{general7eq1}) and (\ref{general8eq1}).
Since (\ref{holdiscserminKBessellemmaeq1}) satisfies (\ref{general6eq1}),
any other solution $c(\zeta,\varphi_1,\varphi_2)$ of (\ref{general6eq1}) is of the form
\begin{align}\label{holdiscserminKBessellemmaeq3}
  \nonumber c(\zeta,\varphi_1,\varphi_2)&=c_3(\varphi_1,\varphi_2)
  \Big(\zeta(\cos(\varphi_1)\cos(\varphi_2)+i\sin(\varphi_1)\sin(\varphi_2))
  +\zeta^{-1}(\cos(\varphi_1)\sin(\varphi_2)+i\sin(\varphi_1)\cos(\varphi_2))\Big)^{m}\\
  &\qquad\Big((\zeta^2-\zeta^{-2})(\cos(2\varphi_1)
  +i\sin(2\varphi_1)\sin(2\varphi_2))+(\zeta^2+\zeta^{-2})\cos(2\varphi_2)\Big)
  ^{\frac{l-l'-m}2}.
\end{align}
with a function $c_3(\varphi_1,\varphi_2)$. Substituting
(\ref{holdiscserminKBessellemmaeq3}) into (\ref{general7eq1}) and
simplifying, we obtain $\frac{\partial c_3}{\partial\varphi_1}=0$.
Substituting (\ref{holdiscserminKBessellemmaeq3}) into
(\ref{general8eq1}) and simplifying, we obtain
$\frac{\partial c_3}{\partial\varphi_2}=0$. Hence $c_3$ is constant. This proves
the uniqueness statement.\qed
\nl
Note that we have shown that if $B$ is a non-zero function
in $\mathcal{S}(\Lambda,\theta,l,l')$,
then, in a neighborhood of the identity, it is given by the unique function obtained from
(\ref{general5eq1}) and Lemma \ref{holdiscserminKBessellemma}.
\subsubsection*{Necessary condition for existence}
The next lemma states
that the analyticity of $B$ puts restrictions on the possible characters $\Lambda$.
\begin{lemma}\label{Lambda-restriction-lemma}
 Let $\Lambda$ be a character of $T(\R)$ defined in (\ref{Lambdaexpliciteq}),
 with $l+l'+m$ even. A necessary condition for $\mathcal{S}(\Lambda,\theta,l,l')$
 to be non-zero is $|m| \leq l-l'$.
\end{lemma}
{\bf Proof:} Assume there exists a non-zero element $B\in\mathcal{S}(\Lambda,\theta,l,l')$.
Being $K$-finite, $B$ is analytic.
As we saw above, in a neighborhood of the identity $B$ is, up to a scalar, given by
\begin{equation}\label{holdiscserminKBesselthmeq1}\renewcommand{\arraystretch}{1.4}
 B(h(\lambda,\zeta,\varphi_1,\varphi_2))=\left\{\begin{array}{l@{\qquad\text{if }}l}
 c_1(\zeta,\varphi_1,\varphi_2)\lambda^{\frac{l+l'+s}2}
 e^{-2\pi\lambda(\zeta^2+\zeta^{-2})}&\lambda>0,\\
 0&\lambda<0,\end{array}\right.
\end{equation}
with $c_1$ as in (\ref{holdiscserminKBessellemmaeq1}) or
(\ref{holdiscserminKBessellemmaeq2}). It follows from (\ref{holdiscserminKBesselthmeq1})
that there exists an analytic function $C_1$ on $\R_{>0}\times(-\pi,\pi)\times(-\pi,\pi)$ given by
$(\zeta,\varphi_1,\varphi_2)\mapsto c_1(\zeta,\varphi_1,\varphi_2)$ on the domain of the function $c_1(\zeta,\varphi_1,\varphi_2)$. Now, first suppose that $m > l-l'$. Then, from (\ref{holdiscserminKBessellemmaeq1}), we see that $(1,0,\pi/4)$ is a limit point of the domain of $c_1(\zeta,\varphi_1,\varphi_2)$, which goes to infinity as the argument approaches the limit point. But that contradicts the analyticity of $C_1$. If $m < -(l-l')$, then we can get a contradiction by a similar argument if we use (\ref{holdiscserminKBessellemmaeq2}). This completes the proof of the lemma. \qed
\subsubsection*{Extending functions from the neighborhood of identity to the whole group}
In the next lemma, we obtain a function on the whole group $\GSp(4,\R)$ that has a
specified behavior in a neighborhood of the identity. This function will be used below
to construct an element $B \in\mathcal{S}(\Lambda,\theta,l,l')$.
\begin{lemma}\label{wfunctionlemma}
 The function $w:\:\GSp(4,\R)\rightarrow\C$ defined by
 \begin{equation}\label{wfunctionlemmaeq3}
 w(h)=i\det(J(h,\mat{i}{}{}{i}))\det(J(h,\mat{-i}{}{}{i}))
  \big({\rm tr}(h\langle \mat{-i}{}{}{i} \rangle)-{\rm tr}
  (h\langle\mat{i}{}{}{i}\rangle)\big)
\end{equation}
satisfies the following properties.
 \begin{enumerate}
  \item $w$ is polynomial in the entries of $h=(h_{ij})\in\GSp(4,\R)$.
  \item If $h=\gamma tuh(\lambda,\zeta,\varphi_1,\varphi_2)r_3(\varphi_3)r_4(\varphi_4)$
   with $\gamma\in\R_{>0}$, $t\in T^1(\R)$, $u\in U(\R)$, $\lambda\in\R^\times$,
   $\zeta\in\R_{>0}$ and $\varphi_1,\ldots,\varphi_4\in\R$, then
   \begin{equation}\label{wfunctionlemmaeq1}
    w(h)=\gamma^4\lambda e^{-2i\varphi_4}
    \Big((\zeta^2-\zeta^{-2})(\cos(2\varphi_1)+i\sin(2\varphi_1)\sin(2\varphi_2))
    +(\zeta^2+\zeta^{-2})\cos(2\varphi_2)\Big).
   \end{equation}
  \item For $\gamma\in\R^\times$, $t\in T^1(\R)$, $u\in U(\R)$,
   $\varphi_3,\varphi_4\in\R$ and $h\in\GSp(4,\R)$ we have
   \begin{equation}\label{wfunctionlemmaeq2}
    w\big(\gamma tuhr_3(\varphi_3)r_4(\varphi_4)\big)=\gamma^4e^{-2i\varphi_4}w(h).
   \end{equation}
 \end{enumerate}
\end{lemma}
{\bf Proof:} i) Note that ${\rm tr}(h\langle\mat{i}{}{}{i}\rangle)$ and
${\rm tr}(h\langle\mat{-i}{}{}{i}\rangle)$ are rational functions in the entries of $h$.
In fact, ${\rm tr}(h\langle\mat{-i}{}{}{i}\rangle)$ has zeros in the denominator
and is not everywhere defined. However, all denominators are cancelled by the
determinant factors, so that we obtain an everywhere defined polynomial function
in the entries of $h$.
\nll
ii) Let $h=\gamma tuh(\lambda,\zeta,\varphi_1,\varphi_2)r_3(\varphi_3)r_4(\varphi_4)$
with $u=\mat{1}{X}{}{1}$. Calculations show that
\begin{itemize}
 \item $\det(J(h,\mat{i}{}{}{i})) = \gamma^2 e^{-i\varphi_3}e^{-i\varphi_4}$,
 \item $\det(J(h,\mat{-i}{}{}{i})) = \gamma^2\cos(2\varphi_2)e^{i\varphi_3}e^{-i\varphi_4}$,
 \item ${\rm tr}(h\langle\mat{i}{}{}{i}\rangle) = i\lambda(\zeta^2+\zeta^{-2})+{\rm tr}(X)$,
 \item ${\rm tr}(h\langle \mat{-i}{}{}{i} \rangle) =\frac{-i\lambda(\zeta^2-
  \zeta^{-2})}{\cos(2\varphi_2)}(\cos(2\varphi_1)+i\sin(2\varphi_1)\sin(2\varphi_2))
  +{\rm tr}(X)$.
\end{itemize}
This proves formula (\ref{wfunctionlemmaeq1}).
\nll
iii) The transformation property (\ref{wfunctionlemmaeq2}) is easily verified from
(\ref{wfunctionlemmaeq3}). Note that the elements $r_3(\varphi_3)$
and $r_4(\varphi_4)$ fix both $\mat{i}{}{}{i}$ and $\mat{-i}{}{}{i}$,
whereas the elements $r_1(\varphi_1)$ and $r_2(\varphi_2)$ only fix $\mat{i}{}{}{i}$.
\qed
\nl
We now state the main theorem about the existence and uniqueness of non-split
Bessel functions.
\begin{theorem}\label{holdiscserminKBesselthm}
 Let $l\geq l'>0$ be integers. Let $\Lambda$ be the
 character of $T(\R)$ defined by a pair $(s,m)\in\C\times\Z$, as in
 (\ref{Lambdaexpliciteq2}).
 If $l+l'+m$ is odd, then $\mathcal{S}(\Lambda,\theta,l,l')=0$. If $l+l'+m$ is even, then
  \begin{equation}\label{dim-condition}
   \dim_\C(\mathcal{S}(\Lambda,\theta,l,l'))=\left\{\begin{array}{l@{\qquad\text{if }}l}
   1& |m| \leq l-l',\\0&|m|>l-l'.\end{array}\right.
  \end{equation}
 Assume that $l+l'+m$ is even and $|m| \leq l-l'$, and let $B_0$ be a non-zero element of
 $\mathcal{S}(\Lambda,\theta,l,l')$. Then, for all $\lambda\in\R^\times$,
 $\zeta>0$ and $\varphi_1,\varphi_2\in\R$, we have, up to a scalar,
 \begin{equation}\label{holdiscserminKBesselthmeq2}\renewcommand{\arraystretch}{1.4}
  B_0(h(\lambda,\zeta,\varphi_1,\varphi_2))=\left\{\begin{array}{l@{\qquad\text{if }}l}
  c_1(\zeta,\varphi_1,\varphi_2)\lambda^{\frac{l+l'+s}2}
  e^{-2\pi\lambda(\zeta^2+\zeta^{-2})}&\lambda>0,\\
  0&\lambda<0,\end{array}\right.
 \end{equation}
 where $c_1$ is the function given in Lemma \ref{holdiscserminKBessellemma}.
 Here, $h(\lambda,\zeta,\varphi_1,\varphi_2)$ is as in (\ref{helementdefeq}).
 Moreover, there exist analytic functions $A_j:\:K^1\rightarrow\C$, for
 $j\in\{-(l-l'),\ldots,l-l'\}$, which are zero for $j\not\equiv l-l'$ mod $2$,
 such that
 \begin{equation}\label{holdiscserminKBesselthmeq3}\renewcommand{\arraystretch}{1.4}
  B_0(h(\lambda,\zeta,0,0)k)=\Big(\sum_{j=-(l-l')}^{l-l'}A_j(k)\zeta^j\Big)
  \lambda^{\frac{l+l'+s}2}e^{-2\pi\lambda(\zeta^2+\zeta^{-2})}
 \end{equation}
 for all $\lambda,\zeta>0$ and $k\in K^1$.
\end{theorem}
{\bf Proof:} We saw that, up to a scalar, any non-zero element $B$ of
$\mathcal{S}(\Lambda,\theta,l,l')$ coincides with
the function (\ref{holdiscserminKBesselthmeq2}) in a neighborhood of the identity; note
that, by (\ref{HRBesseldecompeq}), the transformation properties of $B$ determine it
on elements of the form $h(\lambda,\zeta,\varphi_1,\varphi_2)$.
Since analytic functions on the identity component of $\GSp(4,\R)$ are determined
in a neighborhood of the identity, it follows that
$\dim_\C(\mathcal{S}(\Lambda,\theta,l,l'))\leq1$. It also follows that
(\ref{holdiscserminKBesselthmeq2}) holds for \emph{all} $\lambda\in\R^\times$,
$\zeta>0$ and $\varphi_1,\varphi_2\in\R$. Looking at the formula
(\ref{holdiscserminKBessellemmaeq1}) for $c_1$ and assuming that $|m|\leq l-l'$,
we obtain (\ref{holdiscserminKBesselthmeq3}) in a neighborhood of the identity,
and then by analyticity in general.
\nl
In Lemma \ref{Lambda-restriction-lemma}, we have already shown that $\mathcal{S}(\Lambda,\theta,l,l')=0$ if $|m| > l-l'$.  Assuming that $|m|\leq l-l'$, it remains to find a function satisfying the conditions
($\mathcal{S}1$) -- ($\mathcal{S}5$) defining the space $\mathcal{S}(\Lambda,\theta,l,l')$.
Let $h=(h_{ij})\in\GSp(4,\R)$.
Let $w$ be the polynomial function from Lemma \ref{wfunctionlemma}.
If $0\leq m\leq l-l'$, we define
\begin{equation}\label{bess-model-no-coord-general1}
  B(h)=\mu_2(h)^{l'+\frac s2 + \frac m2}
  \det(J(h,iI_2))^{-l}w(h)^{\frac{l-l'-m}2}\big(h_{44}-h_{32} + ih_{42} + ih_{34}\big)^m
  e^{2 \pi i {\rm tr}(h\langle iI_2 \rangle)}
\end{equation}
if $\mu_2(h) > 0$, and $B(h) = 0$ if $\mu_2(h) < 0$.
If $-(l-l')\leq m\leq0$, we define
\begin{equation}\label{bess-model-no-coord-general2}
  B(h)=\mu_2(h)^{l'+\frac s2 - \frac m2}
  \det(J(h,iI_2))^{-l}w(h)^{\frac{l-l'+m}2}\big(h_{44}-h_{32} + ih_{42} + ih_{34}\big)^{-m}
  e^{2 \pi i {\rm tr}(h\langle iI_2 \rangle)}
\end{equation}
if $\mu_2(h) > 0$, and $B(h) = 0$ if $\mu_2(h) < 0$.
Then $B$ is a well-defined analytic function on $\GSp(4,\R)$, and we shall prove
it satisfies the conditions ($\mathcal{S}1$) -- ($\mathcal{S}5$).
Conditions ($\mathcal{S}2$) and ($\mathcal{S}4$) are verified by a direct calculation,
observing iii) of Lemma \ref{wfunctionlemma}.
Let $h\in\GSp(4,\R)$. In a neighborhood of the identity we can write
$h=tuh(\lambda,\zeta,\varphi_1,\varphi_2)r_3(\varphi_3)r_4(\varphi_4)$
with $t\in T(\R)$, $u\in U(\R)$, $\lambda$ and $\zeta$ close to $1$ and
$\varphi_1,\ldots,\varphi_4$ close to $0$; see (\ref{HRBesseldecompeq}).
Let $t =\mat{g}{}{}{\det(g){}^{t}g^{-1}}$ where $g=\gamma
\mat{\cos(\delta)}{\sin(\delta)}{-\sin(\delta)}{\cos(\delta)}$ and
let $u=\mat{1}{X}{0}{1}$. Calculations show that
\begin{enumerate}
 \item $\mu_2(h) = \gamma^2 \lambda$,
 \item $\det(J(h,iI_2)) = \gamma^2 e^{-i\varphi_3}e^{-i\varphi_4}$,
 \item ${\rm tr}(h\langle iI_2 \rangle) = i\lambda(\zeta^2+\zeta^{-2})+{\rm tr}(X)$,
 \item $h_{44}-h_{32} + ih_{42} + ih_{34} = \gamma e^{i\delta}
  e^{-i\varphi_4}\Big(\zeta(\cos(\varphi_1)\cos(\varphi_2)+i\sin(\varphi_1)\sin(\varphi_2))
  +\zeta^{-1}(\cos(\varphi_1)\sin(\varphi_2)+i\sin(\varphi_1)\cos(\varphi_2))\Big)$,
 \item $h_{44}+h_{32} + ih_{42} - ih_{34} = \gamma e^{-i\delta}
  e^{-i\varphi_4}\Big(\zeta(\cos(\varphi_1)\cos(\varphi_2)+i\sin(\varphi_1)\sin(\varphi_2))
  -\zeta^{-1}(\cos(\varphi_1)\sin(\varphi_2)+i\sin(\varphi_1)\cos(\varphi_2))\Big)$.
\end{enumerate}
Substituting the above formulas into
(\ref{bess-model-no-coord-general1}) (resp.\ (\ref{bess-model-no-coord-general2})),
we conclude that $B$, as defined in (\ref{bess-model-no-coord-general1})
(resp.\ (\ref{bess-model-no-coord-general2})), is given by the formula
(\ref{holdiscserminKBesselthmeq1}) in a neighborhood of the identity.
It follows that $B$ satisfies ($\mathcal{S}5$) in a neighborhood of the identity,
and then, by analyticity, everywhere. The condition $N_+.B=0$ implies that
$B$ is $K$-finite.
The functions $\mu_2(h)$ and $1/\det(J(h,iI_2))$ are slowly increasing.
Since $w(h)$ is a polynomial function,
it follows from (\ref{bess-model-no-coord-general1}) and
(\ref{bess-model-no-coord-general2}) that $B(h)$ is slowly increasing.
This concludes the proof that $B\in\mathcal{S}(\Lambda,\theta,l,l')$.\qed
\subsection{Non-split Bessel models for lowest weight representations}\label{nonsplitmodelssec}
In this section, we will use Theorem \ref{holdiscserminKBesselthm} to obtain the existence and uniqueness of Bessel models for the lowest weight representations of $\GSp(4,\R)$. Let $\mathcal{U}(\mathfrak{g}^1_\C)$ be the universal enveloping algebra of
$\mathfrak{g}^1_\C=\mathfrak{sp}(4,\C)$. If $B$ is any $K^1$-finite function on
$\GSp(4,\R)$, then $\mathcal{U}(\mathfrak{g}^1_\C)B$ is a $(\mathfrak{g}^1, K^1)$-module.
Similarly, if $B$ is any $K$-finite function on
$\GSp(4,\R)$, then $\mathcal{U}(\mathfrak{g}_\C)B$ is a $(\mathfrak{g},K)$-module.
Assume that $B_0$ is the non-zero element of $\mathcal{S}(\Lambda,\theta,l,l')$
described in Theorem \ref{holdiscserminKBesselthm}. Since $\mathfrak{p}_-.B=N_+.B=0$ and
$Z.B$ and $Z'.B$ are multiples of $B$, it follows
that $\mathcal{U}(\mathfrak{g}_\C)B_0$ is spanned by functions of the form
\begin{equation}\label{PBWconsequenceeq}
 X_+^\alpha P_{1+}^\beta P_{0+}^\gamma N_-^\delta B_0,\qquad
 \alpha,\beta,\gamma,\delta\geq0,\;\delta\leq l-l'.
\end{equation}
The main ingredient in the proof of the existence and uniqueness of the Bessel models is the irreducibility of $\mathcal{U}(\mathfrak{g}_\C)B_0$. We will first state a few lemmas which will be used in the proof of the irreducibility of $\mathcal{U}(\mathfrak{g}_\C)B_0$.

\subsection*{General shape of elements in $\mathcal{U}(\mathfrak{g}_\C)B_0$}
\begin{lemma}\label{generalformlemma}
 Let $\Lambda$ be the character of $T(\R)$ defined by a pair $(s,m)\in\C\times\Z$, as in
 (\ref{Lambdaexpliciteq2}). Let $l\geq l'>0$ be integers. We assume that $l+l'+m$ is even
 and that $|m|\leq l-l'$, so that $\mathcal{S}(\Lambda,\theta,l,l')$ is one-dimensional;
 see Theorem \ref{holdiscserminKBesselthm}. Let $B_0$ be a function spanning this space.
 Then every function $B\in\mathcal{U}(\mathfrak{g}_\C)B_0$ is zero on
 the non-identity component of $\GSp(4,\R)$. On the identity component, the
 function $B(h(\lambda,\zeta,0,0)k)$,
 $\lambda,\zeta>0$, $k\in K^1$, is a linear combination of functions of the form
 \begin{equation}\label{generalformlemmaeq1}
  \Big(\sum_{j=-(l-l')}^{l-l'}A_j(k)\zeta^j\Big)(\zeta^2+\zeta^{-2})^e
  (\zeta^2-\zeta^{-2})^fP(\lambda)\lambda^{e+f+\frac{l+l'+s}2}
  e^{-2\pi\lambda(\zeta^2+\zeta^{-2})},
 \end{equation}
 where each $A_j$ is an analytic function on $K^1$, zero if $j\not\equiv l-l'$ mod $2$,
 where $e,f\in\Z$, $e\geq0$, and where $P$ is a polynomial in $\lambda$
 (the functions $A_j$ and $P$, and the exponents $e,f$ all depend on $B$).
\end{lemma}
{\bf Proof:} By (\ref{holdiscserminKBesselthmeq3}), at least the vector $B_0$
has the asserted property (with $e=f=0$).
To prove it in general, it is enough to show that if $B$ has
the asserted property, and if $X$ is one of the operators $X_\pm$, $P_{1\pm}$, $P_{0\pm}$
or $N_\pm$, then $X.B$ also has the asserted property; see (\ref{PBWconsequenceeq}).
It is clear that if $B$ is zero on the non-identity component, then $X.B$ is as well.
To prove that $(X.B)(h(\lambda,\zeta,0,0)k)$ is a linear combination of functions
of the form (\ref{irreducibilitylemmaeq1}), we investigate each term in the formulas
(\ref{N-formula}) -- (\ref{P0-formula}) for the differential operators.
Observing that
\begin{align*}
 \lambda\frac\partial{\partial\lambda}e^{-2\pi\lambda(\zeta^2+\zeta^{-2})}
  &=-2\pi\lambda(\zeta^2+\zeta^{-2})e^{-2\pi\lambda(\zeta^2+\zeta^{-2})},\\
 \zeta\frac\partial{\partial\zeta}\big((\zeta^2+\zeta^{-2})^e(\zeta^2-\zeta^{-2})^f\big)
  &=2e(\zeta^2+\zeta^{-2})^{e-1}(\zeta^2-\zeta^{-2})^{f+1}
  +2f(\zeta^2+\zeta^{-2})^{e+1}(\zeta^2-\zeta^{-2})^{f-1},\\
 \zeta\frac\partial{\partial\zeta}e^{-2\pi\lambda(\zeta^2+\zeta^{-2})}
  &=-4\pi\lambda(\zeta^2-\zeta^{-2})e^{-2\pi\lambda(\zeta^2+\zeta^{-2})},
\end{align*}
and that the operator $\lambda\frac\partial{\partial\lambda}$
(resp.\ $\zeta\frac\partial{\partial\zeta}$) does not change the degree of
a $\lambda$-polynomial (resp.\ $\zeta$-polynomial),
it follows that each time the exponent of one of $\zeta^2+\zeta^{-2}$ or
$\zeta^2-\zeta^{-2}$ is increased, the other is decreased or the exponent of $\lambda$
is increased. This completes the proof.\qed
\nl
{\bf Remark:} The proof shows slightly more, namely that if $B$ is the function
in (\ref{PBWconsequenceeq}), then $B(h(\lambda,\zeta,0,0)k)$ is a linear combination
of functions of the form (\ref{generalformlemmaeq1}), where
$e+f+\deg(P)\leq\alpha+\beta+\gamma$.
\subsubsection*{Convergence of certain integrals}
\begin{lemma}\label{doubleintegrallemma}
 Let $\alpha,\delta\in\R$ and $\beta,\gamma\in\Z$. Then the integral
 \begin{equation}\label{doubleintegrallemmaeq1}
  \int\limits_0^\infty\int\limits_1^\infty\zeta^\alpha
  (\zeta^2-\zeta^{-2})^\beta(\zeta^2+\zeta^{-2})^\gamma\lambda^\delta
  e^{-4\pi\lambda(\zeta^2+\zeta^{-2})}\,d\zeta\,d\lambda
 \end{equation}
 is convergent if and only if
 $\delta>-1$ and $\alpha+2(\beta+\gamma)<2\delta+1$.
\end{lemma}
{\bf Proof:} Since the function is positive, we may consider iterated integrals
in any order. Carrying out the $\lambda$-integration first, we see that this
integral is divergent if $\delta\leq-1$. Assume that $\delta>-1$. Then
$$
  \int\limits_0^\infty\lambda^\delta e^{-4\pi\lambda(\zeta^2+\zeta^{-2})}\,d\lambda
  =(4\pi)^{-(\delta+1)}\Gamma(\delta+1)(\zeta^2+\zeta^{-2})^{-(\delta+1)}.
$$
Hence (\ref{doubleintegrallemmaeq1}) is convergent if and only if
\begin{equation}\label{doubleintegrallemmaeq2}
 \int\limits_1^\infty\zeta^\alpha
  (\zeta^2-\zeta^{-2})^\beta(\zeta^2+\zeta^{-2})^{\gamma-\delta-1}\,d\zeta<\infty.
\end{equation}
Note that $\int_1^2(\zeta^2-\zeta^{-2})^\beta\,d\zeta$ is finite for \emph{all}
integers $\beta$, so that the behavior at $\infty$ determines the convergence.
The integral (\ref{doubleintegrallemmaeq2}) therefore converges if and only if
$\int_1^\infty\zeta^{\alpha+2(\beta+\gamma-\delta-1)}\,d\zeta$ converges.
This is the case if and only if the exponent is less than $-1$. The assertion
follows.\qed
\nl
Next we consider scalar products of Bessel functions. Note that if $\Lambda$ is a
unitary character, and $B_1,B_2$ have the $(\Lambda,\theta)$ Bessel transformation
property, then the function $B_1(g)\overline{B_2(g)}$ is left $R(\R)$ invariant.
For any measurable function on $R(\R)\backslash G(\R)$ we have the integration
formula
\begin{equation}\label{Rintegrationformulaeq}
 \int\limits_{R(\R)\backslash\GSp(4,\R)}f(g)\,dg
 =\int\limits_{K^1}\int\limits_{\R^\times}\int\limits_1^\infty
 f(\begin{bmatrix}\lambda\zeta\\&\lambda\zeta^{-1}\\&&\zeta^{-1}\\&&&\zeta\end{bmatrix}k)
 \frac{\zeta-\zeta^{-3}}{\lambda^4}\,d\zeta\,d\lambda\,dk;
\end{equation}
see (\ref{HRBesseldecompeq}) and \cite{Fu}, (4.7).
\nll
\begin{lemma}\label{scalarproductconvergencelemma}
 Let $\Lambda$ be a unitary character of $T(\R)$.
 For $i=1,2$ let $l_i\geq l_i'>0$ be integers. Let
 $B_{l_i,l_i'}\in\mathcal{S}(\Lambda,\theta,l_i,l_i')$, and let
 $B_i\in\mathcal{U}(\mathfrak{g}_\C)B_{l_i,l_i'}$. If $l_1'+l_2'>4$, then
 \begin{equation}\label{scalarproductconvergencelemmaeq1}
  \int\limits_{R(\R)\backslash\GSp(4,\R)}B_1(g)\overline{B_2(g)}\,dg
 \end{equation}
 is absolutely convergent.
\end{lemma}
{\bf Proof:} By Lemma \ref{generalformlemma}, we may assume that
\begin{equation}\label{scalarproductconvergencelemmaeq2}
 B_i(h(\lambda,\zeta,0,0)k)=
  \Big(\sum_{j=-(l_i-l_i')}^{l_i-l_i'}A_{i,j}(k)\zeta^j\Big)(\zeta^2+\zeta^{-2})^{e_i}
  (\zeta^2-\zeta^{-2})^{f_i}P_i(\lambda)\lambda^{e_i+f_i+\frac{l_i+l_i'+s}2}
  e^{-2\pi\lambda(\zeta^2+\zeta^{-2})},
\end{equation}
where $A_{i,j}$ are analytic functions on $K^1$, and where $P_i$ is a polynomial.
Since $s\in i\R$ and we are taking absolute values, we may assume that $s=0$.
Furthermore, by (\ref{Rintegrationformulaeq}), we may in fact assume that
\begin{equation}\label{scalarproductconvergencelemmaeq3}
 B_i(h(\lambda,\zeta,0,0)k)=
  \zeta^{j_i}(\zeta^2+\zeta^{-2})^{e_i}
  (\zeta^2-\zeta^{-2})^{f_i}P_i(\lambda)\lambda^{e_i+f_i+\frac{l_i+l_i'}2}
  e^{-2\pi\lambda(\zeta^2+\zeta^{-2})},
\end{equation}
where $j_i\leq l_i-l_i'$. The relevant integral is then
$$
 \int\limits_0^\infty\int\limits_1^\infty
  \zeta^{j_1+j_2}(\zeta^2+\zeta^{-2})^{e_1+e_2}
  (\zeta^2-\zeta^{-2})^{f_1+f_2}P_1(\lambda)\overline{P_2(\lambda)}
  \lambda^{e_1+f_1+e_2+f_2+\frac{l_1+l_1'+l_2+l_2'}2}
  e^{-4\pi\lambda(\zeta^2+\zeta^{-2})}
  \frac{\zeta-\zeta^{-3}}{\lambda^4}\,d\zeta\,d\lambda.
$$
By Lemma \ref{doubleintegrallemma}, this integral is finite if and only if
\begin{equation}\label{scalarproductconvergencelemmaeq4}
 e_1+f_1+e_2+f_2+\frac{l_1+l_1'+l_2+l_2'}2-4>-1
\end{equation}
and
\begin{equation}\label{scalarproductconvergencelemmaeq5}
 j_1+j_2-1+2(e_1+e_2+f_1+f_2+1)<2(e_1+f_1+e_2+f_2+\frac{l_1+l_1'+l_2+l_2'}2-4)+1.
\end{equation}
Assuming that $l_1'+l_2'>4$, condition (\ref{scalarproductconvergencelemmaeq4})
is certainly satisfied.
Since $j_i\leq l_i-l_i'$, condition (\ref{scalarproductconvergencelemmaeq5})
is also satisfied provided that $l_1'+l_2'>4$. This concludes the proof.\qed
\subsubsection*{Irreducibility of lowest weight modules}
We will apply Lemma \ref{scalarproductconvergencelemma} to the case where $\mathcal{U}(\mathfrak{g}_\C)B_{l_2,l_2'} \subset \mathcal{U}(\mathfrak{g}_\C)B_{l_1,l_1'}$. In this case, we will show in the following lemma that we always have $l_1'+l_2' > 4$. It will be clear from the proof that only the case $l_1' = 1$ needs to be considered. The lemma is purely algebraic and has nothing to do with the fact
that our functions have the Bessel transformation property.
\begin{lemma}\label{irreducibilitylemma}
 Let $\Lambda$ be the character of $T(\R)$ defined by a pair $(s,m)\in\C\times\Z$, as in
 (\ref{Lambdaexpliciteq2}). Let $l>0$ be an integer. We assume that $l+1+m$ is even
 and that $|m|\leq l-1$, so that $\mathcal{S}(\Lambda,\theta,l,1)$ is one-dimensional;
 see Theorem \ref{holdiscserminKBesselthm}. Let $B_0\in\mathcal{S}(\Lambda,\theta,l,1)$
 be a non-zero element.
 \begin{enumerate}
  \item If $l\geq3$, then the vectors
   \begin{align}\label{irreducibilitylemmaeq1}
     &X_+^{\alpha_1} P_{1+}^2B_0,\qquad X_+^{\alpha_2} P_{0+}B_0,\qquad
      X_+^{\alpha_3} P_{1+}N_-B_0,\qquad X_+^{\alpha_4} N_-^2B_0,\\
     \label{irreducibilitylemmaeq2}&X_+^{\alpha_5} P_{1+}B_0,\qquad X_+^{\alpha_6} N_-B_0,\\
     \label{irreducibilitylemmaeq3}&X_+^{\alpha_7} B_0,
   \end{align}
   where $\alpha_i$ are non-negative integers, are linearly independent.
   No linear combination of these vectors with a fixed weight is annihilated by
   all of $X_-$, $P_{1-}$ and $P_{0-}$.
  \item If $l=2$, then the same statement as in i) holds with
   the last vector in (\ref{irreducibilitylemmaeq1}) omitted.
  \item If $l=1$, then the same statement as in i) holds with
   the last two vectors in (\ref{irreducibilitylemmaeq1}) and
   the last vector in (\ref{irreducibilitylemmaeq2}) omitted.
 \end{enumerate}
\end{lemma}
{\bf Proof:} i) We will repeatedly use the identity
\begin{equation}\label{irreducibilitylemmaeq4}
 X_-X_+^\alpha=X_+^\alpha X_--\alpha X_+^{\alpha-1}Z-\alpha(\alpha-1)X_+^{\alpha-1}
 \qquad(\alpha\geq1)
\end{equation}
in $\mathcal{U}(\mathfrak{g}_\C)$. If any linear combination of the vectors in
(\ref{irreducibilitylemmaeq1}) -- (\ref{irreducibilitylemmaeq3})
is zero, then this linear combination must contain only vectors of the same weight.
Since $X_-B_0=0$, and $\langle X_+,Z,X_-\rangle$ is an $\mathfrak{sl}(2)$-triple,
it is clear (or follows easily by induction using (\ref{irreducibilitylemmaeq4}))
that the vectors $X_+^\alpha B_0$, $\alpha\geq0$, are linearly
independent, and that none of these vectors except $B_0$ is annihilated by $X_-$.
Consider a linear combination of the vectors in (\ref{irreducibilitylemmaeq2}),
\begin{equation}\label{irreducibilitylemmaeq5}
 aX_+^\alpha P_{1+}B_0+bX_+^{\alpha+1}N_-B_0=0.
\end{equation}
If $\alpha=0$, then we get $a=b=0$ by applying $P_{1-}$ to
(\ref{irreducibilitylemmaeq5}). Assume that $\alpha>0$. Applying $X_-$ to
(\ref{irreducibilitylemmaeq5}) and using (\ref{irreducibilitylemmaeq4}), we get
\begin{equation}\label{irreducibilitylemmaeq6}
 -a\alpha(l+\alpha)X_+^{\alpha-1}P_{1+}B_0+
 (a-b(\alpha+1)(l+\alpha-1))X_+^\alpha N_-B_0=0.
\end{equation}
It follows that $a=b=0$ by induction on $\alpha$. This proves the linear independence
of the vectors in (\ref{irreducibilitylemmaeq2}). Assuming linear independence,
the same calculations show that no linear combination of the vectors in
(\ref{irreducibilitylemmaeq2}) is annihilated by both $X_-$ and $P_{1-}$.
Finally we have to consider linear combinations of the vectors in
(\ref{irreducibilitylemmaeq1}). The method is the same as above: Writing down
a linear combination of vectors of the same weight and applying $X_-$, $P_{1-}$ and
$P_{0-}$, we get conditions on the coefficients by induction on $\alpha$.
We omit the details.
\nll
ii) and iii) are proved similarly. Note that $N_-^2B_0=0$ for $l=2$ and
$N_-B_0=0$ for $l=1$.\qed
\begin{proposition}\label{irreducibilityprop}
 Let $\Lambda$ be the character of $T(\R)$ defined by a pair $(s,m)\in\C\times\Z$, as in
 (\ref{Lambdaexpliciteq2}). Let $l\geq l'>0$ be integers. We assume that $l+l'+m$ is even
 and that $|m|\leq l-l'$, so that $\mathcal{S}(\Lambda,\theta,l,l')$ is one-dimensional;
 see Theorem \ref{holdiscserminKBesselthm}. Let $B_0$ be a function spanning this space.
 Then the $(\mathfrak{g},K)$-module $\mathcal{U}(\mathfrak{g}_\C)B_0$ is irreducible.
\end{proposition}
{\bf Proof:} After applying a suitable twist, we may assume that $s=0$
(see the considerations leading up to (\ref{twistedmodelseq}); observe that
$(X.B)\tilde{}=X.\tilde B$ for any of our differential operators $X$).
Let $W\subset\mathcal{U}(\mathfrak{g}_\C)B_0$ be a proper invariant subspace;
we will obtain a contradiction. It is easy to see that $W$ contains a weight
vector $\tilde B_0$ that is annihilated by $X_-$, $P_{1-}$, $P_{0-}$ and $N_+$.
Let $(l_2,l_2')$ be the weight of $\tilde B_0$. Since
$\tilde B_0$ is a $(\Lambda,\theta)$-Bessel function,
$\tilde B_0$ spans the one-dimensional space $\mathcal{S}(\Lambda,\theta,l_2,l_2')$.
We have
$$
 \mathcal{U}(\mathfrak{g}_\C)\tilde B_0\subset W\subsetneq\mathcal{U}(\mathfrak{g}_\C)B_0.
$$
Evidently, $l_2'\geq l'$. Since $X_-X_+^\alpha B_0\neq0$ for any
$\alpha>0$ (see (\ref{irreducibilitylemmaeq4})), we have in fact $l_2'>l'$.
By Lemma \ref{irreducibilitylemma}, if $l'=1$, then $l_2'\geq4$.
Hence, $l'+l_2'>4$ in any case. Therefore, by Lemma \ref{scalarproductconvergencelemma},
$$
 \langle B_1,B_2\rangle:=\int\limits_{R(\R)\backslash\GSp(4,\R)}B_1(g)\overline{B_2(g)}\,dg
$$
is absolutely convergent for all $B_1\in\mathcal{U}(\mathfrak{g}_\C)B_0$ and
$B_2\in\mathcal{U}(\mathfrak{g}_\C)\tilde B_0$. We have $\langle X.B_1,B_2\rangle
=\langle B_1,X.B_2\rangle$ for all $X\in\mathfrak{g}$. Let
$V=\{B_1\in\mathcal{U}(\mathfrak{g}_\C)B_0:\:\langle B_1,B_2\rangle=0\text{ for all }
B_2\in\mathcal{U}(\mathfrak{g}_\C)\tilde B_0\}$. Then $V$ is a proper, invariant
subspace of $\mathcal{U}(\mathfrak{g}_\C)B_0$. Since
$\mathcal{U}(\mathfrak{g}_\C)\tilde B_0$ does not contain the $K^1$-type $(l,l')$,
we have $B_0\in V$. Hence $\mathcal{U}(\mathfrak{g}_\C)B_0\subset V\subsetneq
\mathcal{U}(\mathfrak{g}_\C)B_0$, a contradiction.\qed
\subsubsection*{The main result}
We are now ready to state the main result about non-split Bessel models of lowest weight
representations. Recall that if the character $\Lambda$ of $T(\R)$ is given
by a pair $(s,m)\in\C\times\Z$ as in (\ref{Lambdaexpliciteq2}), then the lowest
weight representation $\mathcal{E}_t(l,l')$ of $\GSp(4,\R)$ can only have a
$(\Lambda,\theta)$-Bessel model if $t=s$; see the end of Sect.\ \ref{besselmodelssec}.
\begin{theorem}\label{besselmodeltheorem}
 Let $\Lambda$ be the character of $T(\R)$ defined by a pair $(s,m)\in\C\times\Z$, as in
 (\ref{Lambdaexpliciteq2}). Let $l\geq l'>0$ be integers. Then the lowest weight
 module $\mathcal{E}_s(l,l')$ of $\GSp(4,\R)$ has a $(\Lambda,\theta)$-Bessel model
 as defined in Sect.\ \ref{besselmodelssec} if and only if $l+l'+m$ is even
 and $|m|\leq l-l'$. If a Bessel model exists, then it is unique.
 The Bessel function $B_0$ representating the highest weight vector in the
 minimal $K$-type of $\mathcal{E}_s(l,l')$ is the function described in
 Theorem \ref{holdiscserminKBesselthm}. The general form of other functions in
 the Bessel model of $\mathcal{E}_s(l,l')$ is described in Lemma \ref{generalformlemma}.
\end{theorem}
{\bf Proof:} Assume that $\pi=\mathcal{E}_s(l,l')$ has a $(\Lambda,\theta)$-Bessel model
$\mathcal{B}_\pi(\Lambda,\theta)$.
Let $B_0\in\mathcal{B}_\pi(\Lambda,\theta)$ be a non-zero vector with the properties
(\ref{typeIannpropeq}). Then $B_0$ is an element of the space
$\mathcal{S}(\Lambda,\theta,l,l')$ defined in Sect.\ \ref{hdsrnonsplitsec}.
It follows from Theorem \ref{holdiscserminKBesselthm} that $l+l'+m$ is even
and $|m|\leq l-l'$ (note that the first condition simply expresses the
compatibility of $\Lambda$ with the central character of $\mathcal{E}(l,l')$).
The one-dimensionality of $\mathcal{S}(\Lambda,\theta,l,l')$ stated in
Theorem \ref{holdiscserminKBesselthm} implies the uniqueness of the space
$\mathcal{B}_\pi(\Lambda,\theta)$.
\nll
Conversely, assume that $l+l'+m$ is even and $|m|\leq l-l'$. Let $B_0$ be a function
spanning the space $\mathcal{S}(\Lambda,\theta,l,l')$; see
Theorem \ref{holdiscserminKBesselthm}. By Proposition \ref{irreducibilityprop}, the
$(\mathfrak{g},K)$-module $\mathcal{U}(\mathfrak{g}_\C)B_0$ is irreducible.
Since it contains a vector with the properties (\ref{typeIannpropeq}), it follows that
$\mathcal{U}(\mathfrak{g}_\C)B_0\cong\mathcal{E}_s(l,l')$.
Hence $\mathcal{U}(\mathfrak{g}_\C)B_0$ provides a $(\Lambda,\theta)$-Bessel model
for $\mathcal{E}_s(l,l')$.\qed
\begin{corollary}\label{besselmodeltheoremcor}
 Let $\Lambda$ be a character of $T^1(\R)\cong S_1$ defined by $m\in\Z$; see
 (\ref{Lambdaexpliciteq2}).
 \begin{enumerate}
  \item Let $l\geq l'>0$ be integers. Then the lowest weight
   module $\mathcal{E}(l,l')$ of $\SSp(4,\R)$ has a $(\Lambda,\theta)$-Bessel model
   as defined in Sect.\ \ref{besselmodelssec} if and only if $l+l'+m$ is even
   and $|m|\leq l-l'$. If a Bessel model exists, then it is unique.
   The Bessel function $B_0$ representating the highest weight vector in the
   minimal $K^1$-type of $\mathcal{E}(l,l')$ is the restriction of the
   function described in Theorem \ref{holdiscserminKBesselthm} to
   $\SSp(4,\R)$. The general form of other functions in
   the Bessel model of $\mathcal{E}(l,l')$ is described in Lemma \ref{generalformlemma}.
  \item Let $l'\leq l<0$ be integers. Then the highest weight
   module $\mathcal{E}(l,l')$ of $\SSp(4,\R)$ does not admit
   a $(\Lambda,\theta)$-Bessel model.
 \end{enumerate}
\end{corollary}
{\bf Proof:} This follows from Theorem \ref{besselmodeltheorem} and the considerations
leading up to Proposition \ref{Sp4GSp4prop}.\qed
\nl
Note that we have throughout assumed that $S=\mat{1}{}{}{1}$ is the matrix
defining the character $\theta$ in (\ref{thetadef2eq}). If we would have chosen
$S=\mat{-1}{}{}{-1}$, then it would be the \emph{highest} weight representations
of $\SSp(4,\R)$ that admit Bessel models; see (\ref{pipiepsilonobviouseq}).
\subsubsection*{Realization in $L^p$-spaces}
As a consequence of our explicit formulas for the Bessel functions, we show
that the Bessel models for most lowest weight representations lie in certain
$L^p$-spaces.
\begin{lemma}\label{convergencelemma}
 Let $\Lambda$ be the character of $T(\R)$ defined by a pair $(s,m)\in\C\times\Z$, as in
 (\ref{Lambdaexpliciteq2}). Assume that $s\in i\R$, so that $\Lambda$ is a unitary
 character. Assume further that $l+l'+m$ is even and that $|m|\leq l-l'$, so that
 $\mathcal{S}(\Lambda,\theta,l,l')$ is a one-dimensional space;
 see Theorem \ref{holdiscserminKBesselthm}. Let $B_0$ be a function spanning this
 space.
 \begin{enumerate}
  \item For any $\varepsilon>0$,
   \begin{equation}\label{convergencelemma02}
    \int\limits_{R(\R)\backslash\GSp(4,\R)}|B_0(x)|^{2+\varepsilon}\,dx<\infty
    \qquad\Longleftrightarrow\qquad l'\geq2.
   \end{equation}
  \item We have
   \begin{equation}\label{convergencelemma01}
    \int\limits_{R(\R)\backslash\GSp(4,\R)}|B_0(x)|^2\,dx<\infty
    \qquad\Longleftrightarrow\qquad l'\geq3.
   \end{equation}
 \end{enumerate}
\end{lemma}
{\bf Proof:} Using the formula for $B_0$ from
Theorem \ref{holdiscserminKBesselthm}, this can be proved similarly as
Lemma \ref{doubleintegrallemma}.\qed
\nl
Assuming that the character $\Lambda$ is unitary and $p>0$, let
$L^p(R(\R)\backslash\GSp(4,\R),\Lambda,\theta)$ be the space of all measurable
functions $B:\:\GSp(4,\R)\rightarrow\C$ with the Bessel transformation property
$B(tug)=\Lambda(t)\theta(u)B(g)$ for $t\in T(\R)$ and $u\in U(\R)$, and such that
\begin{equation}\label{Lpintegraleq}
 \int\limits_{R(\R)\backslash\GSp(4,\R)}|B(g)|^p\,dg<\infty.
\end{equation}
\begin{proposition}\label{L2proposition}
 Let $\Lambda$ be the character of $T(\R)$ defined by a pair $(s,m)\in\C\times\Z$, as in
 (\ref{Lambdaexpliciteq2}). Assume that $s\in i\R$, so that $\Lambda$ is a unitary
 character. Let $l\geq l'>0$ be integers. We assume that $l+l'+m$ is even and that
 $|m|\leq l-l'$, so that the lowest weight $(\mathfrak{g},K)$-module
 $\mathcal{E}_s(l,l')$ possesses a $(\Lambda,\theta)$-Bessel model.
 \begin{enumerate}
  \item Let $\varepsilon>0$ be arbitrary. Then the $(\Lambda,\theta)$-Bessel model
   for $\mathcal{E}_s(l,l')$ lies in
   $L^{2+\varepsilon}(R(\R)\backslash\GSp(4,\R),\Lambda,\theta)$ if and only if $l'\geq2$.
  \item The $(\Lambda,\theta)$-Bessel model for $\mathcal{E}_s(l,l')$ lies in
   $L^2(R(\R)\backslash\GSp(4,\R),\Lambda,\theta)$ if and only if $l'\geq3$.
 \end{enumerate}
\end{proposition}
{\bf Proof:} i) By Lemma \ref{convergencelemma}, the lowest weight vector
$B_0\in\mathcal{S}(\Lambda,\theta,l,l')$ lies in
$L^{2+\varepsilon}(R(\R)\backslash\GSp(4,\R),\Lambda,\theta)$ if and only if
$l'\geq2$. Since the Bessel model of $\mathcal{E}_s(l,l')$ is generated by right translates
of $B_0$ and right translation preserves the $L^{2+\varepsilon}$ norm, we have
$\mathcal{E}_s(l,l')\subset L^{2+\varepsilon}(R(\R)\backslash\GSp(4,\R),\Lambda,\theta)$
if and only if $B_0\in L^{2+\varepsilon}(R(\R)\backslash\GSp(4,\R),\Lambda,\theta)$.
This proves i).
\nll
ii) is proved in the same way.\qed
\nl
{\bf Remark:} This result is plausible, given that $\mathcal{E}(l,l')$ is a
discrete series representation for $l'\geq3$, and a limit of discrete series
representation for $l'=2$.
\section{Split Bessel models}
In this section we investigate the existence and uniqueness of split Bessel models
for the lowest and highest weight representations of $\GSp(4,\R)$ and
$\SSp(4,\R)$. As in the non-split case, we shall work with $\GSp(4,\R)$
and use the discussion preceding
Proposition \ref{Sp4GSp4prop} to obtain results for $\SSp(4,\R)$.
After changing models, as explained in Sect.\ \ref{besselmodelssec},
we may throughout assume that
$$
 S=\mat{1}{}{}{-1}.
$$
\subsection{Double coset decomposition}\label{splitdoublecosetsec}
Again we start by deriving representatives for the double coset space
$R(\R)\backslash\GSp(4,\R)/K^1$, where $R(\R)=T(\R)_0U(\R)$.
Let $S=\mat{1}{}{}{-1}$. In this case,
\begin{equation}\label{TReqsplit}
 T(\R)=\{\mat{x}{y}{y}{x}:\:x,y\in\R,\:x^2-y^2 \neq 0\}.
\end{equation}
Let $t_0=\mat{1}{1}{1}{-1}$ and let $A := \{\mat{a}{0}{0}{b}:\:a,b\in\R^\times\}$
be the split torus in $\GL(2,\R)$. We have $T(\R)\simeq A$ via the map
\begin{equation}\label{split torus isomorphism}
 T(\R)\ni\mat{x}{y}{y}{x}\longmapsto t_0^{-1}
\mat{x}{y}{y}{x} t_0 = \mat{x+y}{0}{0}{x-y} \in A.
\end{equation}
Let $N := \{ \mat{1}{\zeta}{0}{1}:\:\zeta\in\R\}$. The Iwasawa
decomposition for $\GL(2,\R)$ implies
\begin{equation}\label{iwasawa decomp}
\GL(2,\R) = A\cdot N\cdot\SO(2) = t_0A\cdot N\cdot \SO(2)
= t_0At_0^{-1}\cdot t_0N\cdot\SO(2) =T(\R)\cdot t_0N\cdot\SO(2).
\end{equation}
Using this and the Iwasawa decomposition for $\GSp(4,\R)$, we get
\begin{equation}\label{disjoint double cosets split2}
 \GSp(4,\R) = T(\R)U(\R)\{\mat{\lambda t_0\mat{1}{\zeta}{0}{1}}{}{}{\det(t_0)^tt_0^{-1}
 \mat{1}{0}{-\zeta}{1}} : (\lambda, \zeta) \in\R^\times \times \R \}K^1.
\end{equation}
Any signs appearing in $T(\R)\cong\R^\times\times\R^\times$ can be absorbed into
$\lambda$ and $K^1$, so that $T(\R)$ can be replaced by $T(\R)_0$. Hence,
\begin{equation}\label{disjoint double cosets split}
 \GSp(4,\R) = R(\R)\{\mat{\lambda t_0\mat{1}{\zeta}{0}{1}}{}{}{\det(t_0)^tt_0^{-1}
 \mat{1}{0}{-\zeta}{1}} : (\lambda, \zeta) \in\R^\times \times \R \}K^1.
\end{equation}
It can be checked that the double cosets in (\ref{disjoint double cosets split})
are disjoint.
Recalling the coordinates (\ref{Kcoordeq}) in a neighborhood of the identity of $K^1$,
we let
\begin{equation}\label{helementdefeqsplit}
 \hat{h}(\lambda,\zeta,\varphi_1,\varphi_2):=
 \mat{\lambda t_0\mat{1}{\zeta}{0}{1}}{}{}{\det(t_0)^tt_0^{-1}
 \mat{1}{0}{-\zeta}{1}}r_1(\varphi_1)r_2(\varphi_2).
\end{equation}
\subsection{Differential operators}\label{diffopsecsplit}
Following the method of Sect.\ \ref{diffopsecnonsplit}, we will now
derive explicit formulas for the differential operators
given by elements of the complexified Lie algebra $\mathfrak{g}_\C$ on
the functions in a split Bessel model.
Assume that $\mathcal{B}_{\Lambda,\theta}(\pi)$ is a Bessel model for the
$(\mathfrak{g},K)$-module $(\pi,V)$.
For any $B\in\mathcal{B}_{\Lambda,\theta}(\pi)$ we define a function
$f=f_B$ on $\R^\times\times \R \times\R\times\R$ by
\begin{equation}\label{fdefeqsplit}
 f(\lambda,\zeta,\varphi_1,\varphi_2)=B(\hat{h}(\lambda,\zeta,\varphi_1,\varphi_2)).
\end{equation}
It follows from (\ref{disjoint double cosets split}) that if $B$
has a weight $(l,l')$, then $B$ is determined by $f$. If $L$
denotes one of the operators $N_\pm,X_\pm,P_{0\pm},P_{1\pm}$, then
$L.B$ will be determined by the associated function $f_{L.B}$.
We first calculate the action of an element $L$ of the
non-complexified Lie algebra $\mathfrak g$, given by
$$
 (L.B)(\hat{h}(\lambda,\zeta,\varphi_1,\varphi_2))=
 \frac d{dt}\Big|_0B\big(\hat{h}(\lambda,\zeta,\varphi_1,\varphi_2)\exp(tL)\big).
$$
At least for small values of $t$, we can decompose the argument
according to (\ref{disjoint double cosets split}),
\begin{align}\label{exptLiwasawaeqsplit}
 \hat{h}(\lambda,\zeta,\varphi_1,\varphi_2)\exp(tL)
 &=\begin{bmatrix}1&&x(t)&y(t)\\&1&y(t)&z(t)\\&&1\\&&&1\end{bmatrix}
  \mat{g(t)}{}{}{\det(g(t))\,^tg(t)^{-1}}\nonumber\\
 &\hspace{20ex}\hat{h}(\lambda(t),\zeta(t),\varphi_1(t),\varphi_2(t))
  \,r_3(\varphi_3(t))r_4(\varphi_4(t)).
\end{align}
Here, $g(t)\in T(\R)_0$, and $x(t)$ etc.\ are smooth functions in a
neighborhood of $0$ satisfying
\begin{align*}
 &x(0)=y(0)=z(0)=\varphi_3(0)=\varphi_4(0)=0,\\
 &\lambda(0)=\lambda,\quad\zeta(0)=\zeta,
 \quad\varphi_1(0)=\varphi_1,\quad\varphi_2(0)=\varphi_2.
\end{align*}
According to (\ref{split torus isomorphism}), we can write
\begin{equation}\label{gtdecompeqsplit}
 g(t)=t_0\mat{a(t)}{0}{0}{b(t)}
 t_0^{-1}
\end{equation}
with smooth functions $a(t)>0$ and $b(t)>0$ such that $a(0)=1$ and
$b(0)=1$. The character $\Lambda$ of $T(\R)_0$ is of the form
\begin{equation}\label{Lambdaexpliciteqsplit}
 \Lambda(t_0\mat{a}{0}{0}{b}t_0^{-1})=a^{s_1}b^{s_2}, \qquad a,b>0,
\end{equation}
with some $s_1, s_2 \in \C$. It follows that
\begin{align}\label{Loperationgeneraleqsplit}
 &(L.B)(\hat{h}(\lambda,\zeta,\varphi_1,\varphi_2))
  =\frac d{dt}\Big|_0\Big(\theta(\begin{bmatrix}1&&x(t)&y(t)\\&1&y(t)&z(t)\\
  &&1\\&&&1\end{bmatrix})\Lambda(g(t))e^{i(l\varphi_3(t)+l'\varphi_4(t))}
  B\big(\hat{h}(\lambda(t),\zeta(t),\varphi_1(t),\varphi_2(t))\big)\Big)\nonumber\\
 &\;=\frac d{dt}\Big|_0\Big(a(t)^{s_1}b(t)^{s_2}e^{2\pi i(x(t)-z(t))}
  e^{i(l\varphi_3(t)+l'\varphi_4(t))}
  f\big(\lambda(t),\zeta(t),\varphi_1(t),\varphi_2(t)\big)\Big)\nonumber\\
 &\;=\Big(s_1a'(0)+s_2b'(0)+i(l\varphi_3'(0)+l'\varphi_4'(0))
  +2\pi i(x'(0)-z'(0))\Big)f(\lambda,\zeta,\varphi_1,\varphi_2)\nonumber\\
 &\quad+
  \lambda'(0)\frac{\partial f}{\partial\lambda}(\lambda,\zeta,\varphi_1,\varphi_2)
  +\zeta'(0)\frac{\partial f}{\partial\zeta}(\lambda,\zeta,\varphi_1,\varphi_2)
  +\varphi_1'(0)\frac{\partial f}{\partial\varphi_1}(\lambda,\zeta,\varphi_1,\varphi_2)
  +\varphi_2'(0)\frac{\partial f}{\partial\varphi_2}(\lambda,\zeta,\varphi_1,\varphi_2).
\end{align}
Thus we need the derivatives at $0$ of the auxiliary
functions $\lambda, \zeta, \ldots$. To get these, we differentiate
the matrix equation (\ref{exptLiwasawaeqsplit}) and put $t=0$.
This yields sixteen linear equations from which the desired
derivatives can be determined. We will refrain from listing all these
derivatives, and instead just state the resulting formulas for the
action of the complexified Lie algebra.
Let us write $h$ for the element $\hat{h}(\lambda,\zeta,\varphi_1,\varphi_2)$.
\begin{align}
 Z.B&=lB,\label{Z-formulasplit}\\
 Z'.B&=l'B, \label{Z'-formulasplit}\\
 (N_\pm.B)(h)&=\frac i2\tan(2\varphi_2)(l'-l)f(\lambda,\zeta,\varphi_1,\varphi_2)
  +\frac1{2\cos(2\varphi_2)}\frac{\partial f}{\partial\varphi_1}
  (\lambda,\zeta,\varphi_1,\varphi_2)\mp\frac i2\frac{\partial f}{\partial\varphi_2}
  (\lambda,\zeta,\varphi_1,\varphi_2), \label{N-formulasplit}\\
 (X_\pm.B)(h) &=\Big(\frac14(s_1-s_2)(\cos(2\varphi_1)
  \pm i\sin(2\varphi_1)\sin(2\varphi_2))-
  \frac14(s_1+s_2)\cos(2\varphi_2)\nonumber \\
 &\;\pm\frac{l(\cos^4(\varphi_2)+\sin^4(\varphi_2))}{2\cos(2\varphi_2)}
  \mp\frac{l'}4\sin(2\varphi_2)\tan(2\varphi_2) \nonumber \\
 &\;+ 2\pi i \lambda\big((\cos(2\varphi_1)-\zeta\sin(2\varphi_1))\sin(2\varphi_2)
  \pm i(\zeta\cos(2\varphi_1)-\zeta\cos(2\varphi_2)+\sin(2\varphi_1))\big)\Big)
  f(\lambda,\zeta,\varphi_1,\varphi_2)\nonumber\\
 &\;+\frac12 \cos(2\varphi_2)\lambda\frac{\partial f}{\partial\lambda}
  (\lambda,\zeta,\varphi_1,\varphi_2)\\
 &\;+ \frac 12\big(-\zeta \cos(2\varphi_1)-\sin(2\varphi_1) \pm
  i \sin(2\varphi_2)(\cos(2\varphi_1)-\zeta \sin(2\varphi_1))\big)
  \frac{\partial f}{\partial \zeta}(\lambda,\zeta,\varphi_1,\varphi_2) \nonumber \\
 &\;+\big(\frac14\sin(2\varphi_1) \mp \frac i4(\cos(2\varphi_1)\cos(2\varphi_2) - 1)
  \tan(2\varphi_2)\big)
  \frac{\partial f}{\partial \varphi_1}(\lambda,\zeta,\varphi_1,\varphi_2)\\
 &\;+\frac14\sin(2\varphi_2) \frac{\partial f}{\partial \varphi_2}
  (\lambda,\zeta,\varphi_1,\varphi_2),\label{Xpmoperatorsplit}\\
 (P_{1\pm}.B)(h)&=\Big(\frac12(s_1-s_2)\sin(2\varphi_1)\cos(2\varphi_2)
  \mp\frac i2(s_1+s_2)\sin(2\varphi_2)\nonumber \\
  &\quad+\frac i2 \sin(2\varphi_2) (l+l') - 4\pi i\lambda(\pm i
  \cos(2\varphi_2)(\cos(2\varphi_1) - \zeta \sin(2\varphi_1)) -
  \zeta \sin(2\varphi_2))\Big)f(\lambda,\zeta,\varphi_1,\varphi_2) \nonumber \\
  &\quad\pm i \sin(2\varphi_2) \lambda\frac{\partial f}{\partial\lambda}
  (\lambda,\zeta,\varphi_1,\varphi_2) + \cos(2\varphi_2)(\cos(2\varphi_1)
  - \zeta\sin(2\varphi_1)) \frac{\partial f}{\partial \zeta}
   (\lambda,\zeta,\varphi_1,\varphi_2)\nonumber\\
  &\quad-\frac 12\cos(2\varphi_1)\cos(2\varphi_2) \frac{\partial f}{\partial\varphi_1}
   (\lambda,\zeta,\varphi_1,\varphi_2)\mp \frac i2 \cos(2\varphi_2)
   \frac{\partial f}{\partial\varphi_2}
   (\lambda,\zeta,\varphi_1,\varphi_2), \label{P1pmoperatorsplit}\\
 (P_{0\pm}.B)(h)&=\Big(\frac14(s_1-s_2)(-\cos(2\varphi_1)
  \pm i\sin(2\varphi_1)\sin(2\varphi_2))-
  \frac14(s_1+s_2)\cos(2\varphi_2) \nonumber \\
  &\;\mp \frac l4\sin(2\varphi_2) \tan(2\varphi_2) \pm
  \frac{l'(\sin^4(\varphi_2)+\cos^4(\varphi_2))}{2\cos(2\varphi_2)} \nonumber \\
  &\;+ 2\pi i \lambda\big((\cos(2\varphi_1)-\zeta\sin(2\varphi_1))\sin(2\varphi_2)
  \mp i(\zeta\cos(2\varphi_1)+\zeta\cos(2\varphi_2)+\sin(2\varphi_1))\big)\Big)f
   (\lambda,\zeta,\varphi_1,\varphi_2)\nonumber\\
  &\;+ \frac12 \cos(2\varphi_2) \lambda\frac{\partial f}{\partial\lambda}
  (\lambda,\zeta,\varphi_1,\varphi_2)\nonumber\\
  &\; + \frac 12\big(\zeta \cos(2\varphi_1) + \sin(2\varphi_1) \pm
  i(\cos(2\varphi_1)-\zeta \sin(2\varphi_1)) \sin(2\varphi_2)\big)
  \frac{\partial f}{\partial \zeta}(\lambda,\zeta,\varphi_1,\varphi_2) \nonumber \\
  &\;-\big(\frac14\sin(2\varphi_1) \pm \frac i4(1+\cos(2\varphi_1)
  \cos(2\varphi_2)) \tan(2\varphi_2)\big) \frac{\partial f}{\partial\varphi_1}
  (\lambda,\zeta,\varphi_1,\varphi_2)\nonumber\\
  &\; + \frac14\sin(2\varphi_2) \frac{\partial f}{\partial\varphi_2}
  (\lambda,\zeta,\varphi_1,\varphi_2).\label{P0pmoperatorsplit}
\end{align}
\subsection{Non-existence of split Bessel models}\label{hdsrsplitsec}
In this section, we will show that the lowest weight representations of
$\GSp(4,\R)$ do not admit split Bessel models. Let $\Lambda$ be the character
of $T(\R)_0$ defined in (\ref{Lambdaexpliciteqsplit}).
Let $l\geq l'>0$ be integers, and let
$\pi=\mathcal{E}(l,l')$ be a lowest weight representation of
$\GSp(4,\R)$ as defined in Sects.\ \ref{lowestweightsec} and \ref{SpGSpsec}.
As in Sect.\ \ref{hdsrnonsplitsec}, let $\mathcal{S}(\Lambda,\theta,l,l')$
be the space of all functions $B:\:\GSp(4,\R)\rightarrow\C$ satisfying
the following conditions.
\begin{enumerate}
 \item[($\mathcal{S}1$)] $B$ is smooth and $K$-finite.
 \item[($\mathcal{S}2$)] $B(tug)=\Lambda(t)\theta(u)B(g)$ for all
  $t\in T(\R)_0$, $u\in U(\R)$, $g\in\GSp(4,\R)$.
 \item[($\mathcal{S}3$)] $B$ is slowly increasing.
 \item[($\mathcal{S}4$)] $Z.B=lB$ and $Z'.B = l'B$. Equivalently,
  $B(gr_3(\varphi_3)r_4(\varphi_4))=e^{i(l\varphi_3+l'\varphi_4)}B(g)$
  for all $\varphi_3,\varphi_4\in\R$, $g\in\GSp(4,\R)$.
 \item[($\mathcal{S}5$)] $N_+.B=X_-.B=P_{1-}.B=P_{0-}.B=0$.
\end{enumerate}
If $B$ is a highest weight vector in the minimal $K$-type in a
Bessel model for $\pi$ of type $(\Lambda,\theta)$, then $B$ satisfies
($\mathcal{S}1$) -- ($\mathcal{S}5$). Given such a $B$, we define the
associated function $f$ as in (\ref{fdefeqsplit}). Note that $f$ is an analytic
function on $\R^\times\times\R\times\R\times\R$. Calculations using formulas
(\ref{N-formulasplit}) -- (\ref{P0pmoperatorsplit}) show that condition ($\mathcal{S}5$)
is equivalent to the following system of differential equations,
\begin{align}
 \label{flambdasplit}\frac{\partial f}{\partial \lambda}&=
  \Big(\frac{l+l'+s_1+s_2}{2\lambda} + 4 \pi \zeta\Big) f,\\
 \label{fzetasplit}\frac{\partial f}{\partial \zeta}&=
  \Big(4 \pi \lambda - \frac{(l-l'+s_1-s_2)(\cos(\varphi_2) \sin(\varphi_1)
  - i \cos(\varphi_1) \sin(\varphi_2))}
  {2\cos(\varphi_2)(\cos(\varphi_1) - \zeta \sin(\varphi_1))+2i(\zeta
  \cos(\varphi_1)+\sin(\varphi_1))\sin(\varphi_2)}\Big)f,\\
 \label{fphi1split}\frac{\partial f}{\partial \varphi_1}&=\Big(
  \frac{(s_1-s_2)\cos(2\varphi_2)+(l-l')(\zeta\sin(2\varphi_1)-\cos(2\varphi_1)
  -i\sin(2\varphi_2)(\sin(2\varphi_1)+\zeta\cos(2\varphi_1))}
  {\zeta\cos(2\varphi_2)+\cos(2\varphi_1)
  (-\zeta+i\sin(2\varphi_2))-\sin(2\varphi_1)(1+i\zeta\sin(2\varphi_2))}\Big)f,\\
 \label{fphi2split}\frac{\partial f}{\partial \varphi_2}&=\Big(
  \frac{-i(s_1-s_2)+(l-l')(-\zeta\sin(2\varphi_2)+i\cos(2\varphi_2)
  (\cos(2\varphi_1)-\zeta\sin(2\varphi_1)))}{\zeta\cos(2\varphi_2)+\cos(2\varphi_1)
  (-\zeta+i\sin(2\varphi_2))-\sin(2\varphi_1)(1+i\zeta\sin(2\varphi_2))}\Big)f.
\end{align}
From (\ref{flambdasplit}) we get
\begin{equation}\label{lambdadiffeqsolnsplit}
 f(\lambda, \zeta, \varphi_1, \varphi_2) = \left\{\begin{array}{l@{\qquad\mbox{if }}l}
 c_1(\zeta, \varphi_1, \varphi_2)\lambda^{\frac{l+l'+s_1+s_2}2}e^{4\pi \lambda \zeta}
 &\lambda>0,\\
 c_2(\zeta, \varphi_1, \varphi_2)(-\lambda)^{\frac{l+l'+s_1+s_2}2}e^{4\pi \lambda \zeta}
 &\lambda<0,\end{array}\right.
\end{equation}
with certain functions $c_1(\zeta, \varphi_1, \varphi_2)$ and
$c_2(\zeta, \varphi_1, \varphi_2)$. Note that $c_1$ and $c_2$ are analytic functions
on all of $\R\times\R\times\R$. Assume that $c_1$ is not constantly zero. Then, by
analyticity, there exists a choice of
$(\zeta,\varphi_1,\varphi_2)\in\R_{>0}\times\R\times\R$ such that
$c_1(\zeta,\varphi_1,\varphi_2)\neq0$. But a look at (\ref{lambdadiffeqsolnsplit})
shows that this would violate the moderate growth condition for $B$; see
(\ref{slowlyincreasingeq2}). Similarly, the assumption that $c_2$ is not
constantly zero also violates moderate growth. This proves the following.
\begin{theorem}\label{holdiscserminKBesselsplitthm}
 Let $S\in M(2\times2,\R)$ be a non-degenerate symmetric matrix with $\det(S)<0$.
 Let $\theta$ be the corresponding character of $U(\R)$ as in (\ref{thetadef2eq}),
 and let $T(\R)$ be the group defined in (\ref{TRthetaeq1}). Let $l\geq l'>0$ be integers.
 Then, for any character $\Lambda$ of $T(\R)_0\cong\R_{>0}\times\R_{>0}$, the space
 $\mathcal{S}(\Lambda,\theta,l,l')$ is zero.
\end{theorem}
\begin{corollary}\label{holdiscserminKBesselsplitthmcor}
 The lowest weight representations $\mathcal{E}(l,l')$ of $\GSp(4,\R)$ do not admit split
 Bessel models. The lowest and highest weight representations
 of $\SSp(4,\R)$ do also not admit split Bessel models.
\end{corollary}
{\bf Proof:} The highest weight vector in the minimal $K$-type of a split Bessel model
for $\mathcal{E}(l,l')$ would be a non-zero element of
$\mathcal{S}(\Lambda,\theta,l,l')$. The assertion about $\SSp(4,\R)$
follows from Proposition \ref{Sp4GSp4prop}.\qed
\nl
{\bf Remark:} The system (\ref{flambdasplit}) -- (\ref{fphi2split}) can be solved
formally. Restricting to the connected component of the domain of $f$ given by
$\lambda>0$, the one-dimensional solution space is spanned by the function
\begin{align}\label{final split coord formula}
 f(\lambda, \zeta, \varphi_1, \varphi_2)&=
  \Big(\cos(\varphi_2)(\cos(\varphi_1) - \zeta \sin(\varphi_1))
  +i(\zeta\cos(\varphi_1)+\sin(\varphi_1))\sin(\varphi_2)\Big)^{\frac{l-l'+s_1-s_2}2}
  \nonumber \\
 &\qquad\times\Big(\sin(\varphi_1)\cos(\varphi_2)
  - i \cos(\varphi_1) \sin(\varphi_2)\Big)^{\frac{l-l'-s_1+s_2}2}
   \lambda^{\frac{l+l'+s_1+s_2}2}e^{4\pi \lambda \zeta}.
\end{align}
\section{An application}\label{applicationsec}
In the previous sections we obtained the formula for the
highest weight vector in the minimal $K$-type of a lowest weight representation.
One of the main uses for such a formula is for explicit computations involving the
archimedean Bessel models. For example, if $F$ is a scalar valued
Siegel modular form and $f$ is a Maa{\ss} form, then the formula
for the archimedean Bessel function (already obtained in
\cite{Su}) was used in \cite{Fu}, \cite{P-S1} and \cite{P-S2} to
obtain an integral representation of the degree-$8$ $L$-function
$L(s, F \times f)$.
\nl
Since we now have the formula for the Bessel function for any lowest weight representation,
and in particular any holomorphic discrete series representation,
we will use it to obtain an integral
representation for $L(s, {\bf F} \times f)$, where ${\bf F}$ is a
vector valued holomorphic Siegel modular form. The vector entering the archimedean
zeta integral will actually \emph{not} directly correspond to the modular form
$\mathbf F$, but will be a vector spanning a certain one-dimensional $K$-type
in the lowest weight representation generated by $\mathbf F$.
We will give an algorithm to obtain the formula
for such a vector and explicitly compute it in some
cases. Then we will briefly recall vector valued Siegel modular
forms and the Bessel models associated with them. Finally, we will
consider the integral representation of Furusawa and compute the
archimedean integral in the vector valued holomorphic Siegel
modular forms case.
\subsection{Finding good vectors}\label{finding-good-vector-sec}
Let $l\geq l'>0$ be integers of the same parity, and consider the lowest weight
representation $\mathcal{E}_s(l,l')$ of $\GSp(4,\R)$.
In earlier sections, we have obtained a formula for the highest weight vector in the minimal
$K$-type for such representations. In this section, we give an algorithm (and some examples)
to find the formula for the vector in the $1$-dimensional $K$-type $(l,l)$.
Let $B_0$ be the lowest weight vector in the $(\Lambda,\theta)$-Bessel model
of $\mathcal{E}_s(l,l')$ as described in Theorem \ref{holdiscserminKBesselthm}.
We denote this function also by $B_{l,l'}$, and define recursively
for $k=1,\ldots,(l-l')/2$
\begin{equation}\label{highest-weight-step-by-step}
 B_{l,l'+2k}=\Big(P_{0+}+\frac1\alpha N_-N_+P_{0+}+\frac1{2\alpha(\alpha+1)}
 N_-^2N_+^2P_{0+}\Big)B_{l,l'+2k-2},\qquad\alpha=l-l'-2k+2.
\end{equation}
Calculations using the multiplication table given in Sect.\ \ref{grpsec} show that
$N_+B_{l,l'+2k}=0$ for all $k$ (use that $N_+^3P_{0+}v = 0$ for highest weight
vectors $v$). Hence, $B_{l,l'+2k}$ is the highest weight vector in the
$K$-type $(l,l'+2k)$.
Let us now assume that $S=\mat{1}{}{}{1}$ and $\Lambda \equiv 1$.
Using the formula for $B_{l,l'}$ from Theorem \ref{holdiscserminKBesselthm} and (\ref{highest-weight-step-by-step}) we
will give the formula for $B_{l,l}$ in a few cases. Note that,
since $B_{l,l}$ lies in a $1$-dimensional $K$-type, it is
completely determined by its values on $h(\lambda,\zeta,0,0)$ (see
(\ref{helementdefeq})). Set $x=(\zeta^2+\zeta^{-2})/2$. For
$\lambda < 0$ we have $B_{l,l}(h(\lambda,\zeta,0,0)) = 0$. For $\lambda> 0$
we have the following formulas.
\begin{description}
\item{$l'=l-2$:}
 \begin{equation}\label{l-(l-2)-one-dim-vector}
  B_{l,l}(h(\lambda,\zeta,0,0)) = 4 e^{-4 \pi \lambda x}
  \lambda^{l-1}\Big(2(l-3)x+8 \pi \lambda\Big)
 \end{equation}
\item{$l'=l-4$:}
 \begin{equation}\label{l-(l-4)-one-dim-vector}
  B_{l,l}(h(\lambda,\zeta,0,0)) = \frac{4}{15} e^{-4 \pi \lambda x}
  \lambda^{l-2}\Big(12(l-4)(l-5)x^2 - 8(l-4)x(8\pi \lambda)+
  2(8\pi \lambda)^2 - 4(l-4)(l-5)\Big).
 \end{equation}
\item{$l'=l-6$:}
 \begin{align}\label{l-(l-6)-one-dim-vector}
  &B_{l,l}(h(\lambda,\zeta,0,0))=\frac{8}{35} e^{-4 \pi \lambda x}
   \lambda^{l-3}\Big(40(l-5)(l-6)(l-7)x^3 - 36(l-5)(l-6)(8\pi \lambda) x^2 \nonumber \\
  &\qquad +12(l-5)(8\pi \lambda)^2 x
   -2(8\pi \lambda)^3 - 24(l-5)(l-6)(l-7)x+12(l-5)(l-6)(8\pi \lambda)\Big).
 \end{align}
\item{$l'=l-8$:}
 \begin{align}\label{l-(l-8)-one-dim-vector}
  &B_{l,l}(h(\lambda,\zeta,0,0))=\frac{16}{315} e^{-4 \pi \lambda x}\lambda^{l-4}
   \Big(560(l-6)(l-7)(l-8)(l-9)x^4- 640(l-6)(l-7)(l-8)(8\pi \lambda) x^3 \nonumber \\
  &\quad+ 288(l-6)(l-7)(8\pi \lambda)^2x^2 - 64(l-6)(8\pi \lambda)^3 x + 8(8\pi \lambda)^4
   - 480 (l-6)(l-7)(l-8)(l-9) x^2 \nonumber \\
  &\quad + 384(l-6)(l-7)(l-8)(8\pi \lambda) x - 96 (l-6)(l-7)(8\pi \lambda)^2
   + 48(l-6)(l-7)(l-8)(l-9)\Big).
 \end{align}
\end{description}
\subsection{Vector valued Siegel modular forms and global Bessel models}\label{vvsec}
Let $\SH_2 := \{Z \in M_2(\C) : \,^tZ = Z,\;{\rm Im}(Z) > 0\}$
be the Siegel upper half plane of degree $2$. Let $\Gamma_2 =\SSp(4,\Z)$.
Let $n$ be an odd, positive integer.
Let $(\rho_0, V)$ be the polynomial (holomorphic), irreducible,
$n$-dimensional representation of $\GL(2,\C)$ for which the center acts trivially.
For a positive integer $l\geq n$, let us denote by $\rho$ the representation of
$\GL(2,\C)$ on $V$ given by $g \mapsto \det(g)^l \rho_0(g)$. A
vector valued Siegel modular form of type $\rho$ is defined as a
holomorphic function ${\bf F} : \SH_2 \rightarrow V$ satisfying
\begin{equation}\label{vect-sieg-defn}
{\bf F}(\gamma \langle Z \rangle) = \rho(C Z + D) {\bf F}(Z),\;\mbox{ where } \gamma = \mat{A}{B}{C}{D} \in \Gamma_2,\;Z \in \SH_2,\;\gamma \langle Z \rangle := (AZ+B)(CZ+D)^{-1}.
\end{equation}
Such a function has a Fourier expansion of the form
\begin{equation}\label{vect-sieg-four-exp}
{\bf F}(Z) = \sum\limits_{S \geq 0} {\bf A}(S) e^{2 \pi i \tr(SZ)},
\end{equation}
where $S$ runs through all semi-integral, semi-positive definite,
symmetric $2 \times 2$ matrices. We say that ${\bf F}$ is a cusp
form if ${\bf A}(S) \neq 0$ only if $S > 0$. We denote the space
of vector valued Siegel cusp forms of type $\rho$ with respect to
$\Gamma_2$ by $S_{\rho}(\Gamma_2)$. Let us assume that ${\bf F}
\in S_{\rho}(\Gamma_2)$ is a Hecke eigenform. We will now
construct the automorphic representation of $\GSp(4,\A)$
corresponding to ${\bf F}$. For $g = g_{\Q} g_{\infty} k_0$, with
$g_{\Q} \in \GSp(4,\Q)$, $g_{\infty} = \mat{A}{B}{C}{D} \in
\GSp(4,\R)^+$ and $k_0 \in \prod\limits_{p < \infty}\GSp(4,\Z_p)$, define
\begin{equation}\label{global-vect-func-defn}
{\bf \Phi}(g) := \mu_2(g_{\infty})^{l+(n-1)/2} \rho(C I + D)^{-1} {\bf F}(g_{\infty} \langle I \rangle),
\end{equation}
where $I = \mat{i}{0}{0}{i}$. Choose a fixed non-zero linear functional $\Psi$ on $V$,
and set
\begin{equation}\label{scalar-val-func}
\Phi(g) := \Psi({\bf \Phi}(g)), \qquad g \in \GSp(4,\A).
\end{equation}
Let $(\pi_{{\bf F}}, V_{{\bf F}})$ be an irreducible subspace of
the $\GSp(4,\A)$-space obtained from right translates of $\Phi$.
Then $\pi_{{\bf F}}$ is an irreducible, cuspidal, automorphic
representation of $\GSp(4,\A)$. Note that $\pi_{\bf F}$ does not
depend on the choice of $\Psi$. If $\pi_{\bf F} = \otimes' \pi_p$,
then, for $p < \infty$, $\pi_p$ is an unramified representation of
$\GSp(4,\Q_p)$, and $\pi_{\infty}$ is the lowest weight representation
$\mathcal{E}(l, l-(n-1))$ (which is a holomorphic discrete series
representation if $l\geq n+2$).
\nl
Let $S$ be a positive definite, semi-integral, symmetric $2 \times
2$ matrix. Let the discriminant of $S$ be given by $d(S) =
-\det(2S) = -D$ and $L = \Q(\sqrt{-D})$. Let $T(\A) \simeq
\A_L^{\times}$ be the adelic points of the group defined in
(\ref{TRthetaeq1}). Let $R(\A) = T(\A) U(\A)$ be the Bessel
subgroup of $\GSp(4,\A)$. Let $\Lambda$ be an ideal class
character of $L$, i.e., a character of
$$
 T(\A)/T(\Q)T(\R) \prod\limits_{p<\infty}(T(\Q_p)\cap \GL(2,\Z_p)).
$$
Let $\psi$ be a character of $\Q\backslash \A$ that is trivial
on $\Z_p$ for all primes $p$ and satisfies $\psi(x)=e^{-2\pi ix}$
for all $x\in\R$.

We define the global Bessel function of type $(\bar\Lambda, -S, \psi)$
associated to $\bar{\phi} \in V_{\bf F}$ by
\begin{equation}\label{global Bessel model defn}
B_{\bar{\phi}}(g) = \int\limits_{Z(\A)R(\Q)\backslash
R(\A)}(\Lambda \otimes \theta)(r)^{-1}\bar{\phi}(rg)dr,
\end{equation}
where $\theta(\mat{1}{X}{}{1})=\psi(\tr(SX))$, $Z(\A)$ is the
center of $\GSp(4,\A)$ and $\bar{\phi}(h) = \overline{\phi(h)}$.
Note that $\pi_{\bf F}$ has a global Bessel model of type
$(\bar\Lambda, -S, \psi)$, or equivalently, a Bessel model of type
$(\bar\Lambda, S, \psi^{-1})$, if there is a $\phi \in V_{\bf F}$ such that
$B_{\bar{\phi}} \neq 0$. The archimedean component of the character $\psi^{-1}$
coincides with the character we fixed in Sect.\ \ref{besselsubgroupssec}, so
that our local theory applies without changes.
Note also that, since $\pi_{\bf F}$ is
irreducible, if $B_{\bar{\phi}} \neq 0$ for some $\phi \in V_{\bf
F}$ then the same is true for all elements of $V_{\bf F}$. We now
make the following important assumptions.
\begin{description}
\item[\bf Assumption $1$:] $\pi_{\bf F}$ has a global Bessel model
of type $(\bar\Lambda, S, \psi^{-1})$ such that $d(S)=-D$ is the fundamental
discriminant of $\Q(\sqrt{-D})$.

\item[Assumption $2$:] $l$ is a multiple of $w(-D)$, the
  number of roots of unity in $\Q(\sqrt{-D})$.
\end{description}

{\bf Remark:} Assumption $1$ implies that for any $S'$ in any
$\SL(2,\Z)$ equivalence class of primitive semi-integral, positive
definite matrices with $d(S') = -D$, we can find a $\Lambda'$ such
that $\pi_{\bf F}$ has a $(\Lambda', S', \psi^{-1})$ global Bessel
model. This can be explained as follows. In \cite[1-26]{Su},
Sugano has obtained the following formula for the vector valued
function ${\bf \Phi}$ defined in (\ref{global-vect-func-defn}),
\begin{align}\label{sugano-vect-bessel-formula}
 B_{\bar{{\bf \Phi}}}(g_{\infty}) &=\int\limits_{Z(\A)R(\Q)\backslash
  R(\A)}(\Lambda \otimes \theta)(r)^{-1}\bar{{\bf \Phi}}(rg_{\infty})dr \nonumber \\
 &=\mu_2(g_{\infty})^{l+(n-1)/2} \overline{\rho(C I + D)^{-1} e^{2 \pi i
  \tr(S(g_{\infty} \langle I \rangle))}  \pi_{\bar{\Lambda}}\Big(\frac 1{h(-D)}
  \sum\limits_{i=1}^{h(-D)} \bar{\Lambda}^{-1}(u_i) {\bf A}(S_i) \Big)},
\end{align}
where $g_{\infty} = \mat{A}{B}{C}{D} \in \GSp(4,\R)^+$,
$\pi_{\Lambda} = \int\limits_{T^1(\R)} \rho(\zeta)^{-1}
\Lambda_{\infty}^{-1}(\zeta) d^{\times} \zeta \in {\rm End}(V)$,
$h(-D)$ is the class number of $\Q(\sqrt{-D})$, the elements
$\{u_i\}$ are representatives of the idele class group of
$\Q(\sqrt{-D})$, and $\{ S_i \}$ is the orbit of $S$ under the
action of the $u_i$'s. It is clear that $B_{\Psi({\bf \bar\Phi})} =
\Psi (B_{\bf \bar\Phi})$, and hence, one can find a $\Psi$ such that
$\pi_{\bf F}$ has a $(\bar\Lambda, S, \psi^{-1})$-Bessel model if and only
if $B_{\bf \bar\Phi} \neq 0$. Now, if $-D$ is the fundamental
discriminant of $\Q(\sqrt{-D})$ then the set $\{S_i\}$ in
(\ref{sugano-vect-bessel-formula}) runs through all the
$\SL_2(\Z)$ equivalence classes of primitive semi-integral,
positive definite matrices with discriminant $-D$. This implies
that the non-vanishing of $B_{\bf \bar\Phi}$ depends only on $-D$ and
$\Lambda$ and not on the specific matrix $S$. Finally, note that
(\ref{sugano-vect-bessel-formula}) implies, in the scalar valued
case, that Assumption $1$ is equivalent to assuming that $F$ has a
non-zero Fourier coefficient $A(S)$, where $S$ satisfies the
condition from Assumption $1$. \nl From the above remark, it makes
sense to consider
\begin{equation}\label{special four coeff matrix}\renewcommand{\arraystretch}{1.2}
 S(-D) = \left\{\begin{array}{l@{\qquad}l}
    \mat{\frac D4}{0}{0}{1}& \hbox{ if } D \equiv 0 \pmod{4}, \\[3ex]
 \mat{\frac{1+D}{4}}{\frac 12}{\frac 12}{1}& \hbox{ if } D \equiv 3 \pmod{4}.
 \end{array}\right.
\end{equation}
\subsection{An integral representation}\label{intrepsec}
Let $N = \prod p^{n_p}$ be a positive integer. We denote the space of Maa{\ss} cusp
forms of weight $l_1 \in \Z$ with respect to $\Gamma_0(N)$ by
$S_{l_1}^M(N)$. A function $f \in S_{l_1}^M(N)$ has the Fourier expansion
\begin{equation}\label{Maass form four exp}
 f(x+iy) = \sum\limits_{n \neq 0} a_n W_{{\rm sgn}(n)\frac{l_1}2,
\frac{ir}2}(4 \pi|n|y)e^{2 \pi i nx},
\end{equation}
where $W_{\nu,\mu}$ is a classical Whittaker function and
$(\Delta_{l_1}+\lambda) f=0$ with $\lambda = 1/4 + (r/2)^2$. Here
$\Delta_{l_1}$ is the Laplace operator defined in Sect.\ 2.3 of
\cite{P}. Let $f \in S_{l_1}^M(N)$ be a  Hecke eigenform. If $ir/2
= (l_2-1)/2$ for some integer $l_2 > 0$, then the cuspidal,
automorphic representation of $\GL(2,\A)$ constructed below is
holomorphic at infinity of lowest weight $l_2$. In this case we
make the additional assumptions that $l_2\leq l$ and $l_2\leq
l_1$, where $l$ is coming from the Siegel cusp form ${\bf F}$ as
in the previous section. Starting from a Hecke eigenform $f$, we
obtain another Maa{\ss} form $f_l \in S_{l}^M(N)$ by applying the
raising and lowering operators as in Sect.\ 2.3 of \cite{P}.
Define a function $\hat{f}$ on $\GL(2,\A)$ by
\begin{equation}\label{maass form lift to group}
 \hat{f}(\gamma_0 m k_0) = \Big(\frac{\gamma i +\delta}{|\gamma i +\delta|}\Big)^{-l}
 f_l\Big(\frac{\alpha i + \beta}{\gamma i + \delta}\Big),
\end{equation}
where $\gamma_0 \in \GL(2,\Q)$, $m =
\mat{\alpha}{\beta}{\gamma}{\delta} \in \GL(2,\R)^+$, $k_0 \in
\prod\limits_{p | N}K^{(1)}(\p^{n_p}) \prod\limits_{p \nmid
N}\GL(2,\Z_p)$. Here, for $p | N$ we have $K^{(1)}(\p^{n_p}) =
\GL(2,\Q_p) \cap\mat{\Z_p^{\times}}{\Z_p}{\p^{n_p}}{\Z_p^\times}$
with $\p = p\Z_p$. Let $\tau_f \cong \otimes'_p \tau_{p}$ be the
irreducible, cuspidal, automorphic representation of $\GL(2,\A)$
generated by $\hat f$. \vskip 0.2in

Let the unitary group $\GU(2,2;L)(\A)$ and its subgroups $P$,
$M^{(1)}$, $M^{(2)}$ and $N$ be as defined in p.\ 190, 192 of
\cite{Fu} (or Sect.\ 2.1 of \cite{P-S1}). As in p.\ 210 of
\cite{Fu} (or Sect.\ 5.2 of \cite{P-S1}), define an Eisenstein
series on $\GU(2,2;L)(\A)$ by
\begin{equation}\label{eisenstein series definition}
 E_{\Lambda}(g,s) = \sum\limits_{\gamma \in P(\Q) \backslash
 \GU(2,2;L)(\Q)}f_{\Lambda}(\gamma g,s)
\end{equation}
where $f_\Lambda$ is as defined in p.\ 209 of \cite{Fu} from
$\hat{f}$. For any vector $\phi \in V_{\bf F}$, we consider the
global integral
\begin{equation}\label{global integral calculation}
 Z(s,\Lambda) = \int\limits_{Z(\A)\GSp(4,\Q)\backslash \GSp(4,\A)}E_{\Lambda}(h,s)\bar{\phi}(h)dh.
\end{equation}
In Theorem 2.4 of \cite{Fu}, the following basic identity has been proved,
\begin{equation}\label{basic identity}
Z(s,\Lambda) = \int\limits_{R(\A) \backslash \GSp(4,\A)}W_{\Lambda}(\eta h,
s)B_{\bar{\phi}}(h) dh,
\end{equation}
where $W_{\Lambda}(g,s)$ and $\eta$ are as defined in p.\ 196 of
\cite{Fu} and $B_{\bar{\phi}}$ is as defined in (\ref{global
Bessel model defn}). Let $\pi_{\bf F}$ be given by its
$(\bar\Lambda, S(-D), \psi^{-1})$-Bessel model. Let $\phi = \otimes \phi_p
\in V_{\bf F}$ be such that for $p < \infty$ we choose $\phi_p$ to
be the unique spherical vector in $\pi_p$ satisfying $\phi_p(1) =
1$ and $\phi_{\infty}$ to be a vector in the $1$-dimensional
$K$-type $(l,l)$ in the lowest weight representation
$\mathcal{E}(l, l-(n-1))$. Recall that we have assumed $n$ to be
odd, so that the $1$-dimensional $K$-type $(l,l)$ really occurs in
$\mathcal{E}(l, l-(n-1))$. Since $\varphi$ is a pure tensor, the
global zeta integral factors,
$$
 Z(s,\Lambda) =\prod\limits_p Z_p(s,\Lambda) = \prod\limits_p
  \int\limits_{R(\Q_p) \backslash \GSp(4,\Q_p)}W_p(\eta h,s)
  B_{\bar{\phi}_p}(h) dh.
$$
For $p < \infty$ the integral $Z_p(s,\Lambda)$ has been evaluated
in \cite{Fu}, \cite{P-S1} and \cite{P-S2}. We will now compute
$Z_{\infty}(s, \Lambda)$. From Sect.\ 4.7 of \cite{Fu}, we have
\begin{align}\label{archintegral1eq}
 Z_{\infty}(s) &= \pi\int\limits_{\R^{\times}}\int\limits_1^{\infty}
  W_{\infty}\Big(\eta \begin{bmatrix}\lambda t_0\mat{\zeta}{}{}{\zeta^{-1}}&\\
  &^tt_0^{-1}\mat{\zeta^{-1}}{}{}{\zeta}\end{bmatrix},s\Big) \nonumber \\
 &\hspace{10ex}B_{\bar{\phi}_\infty}\Big(\begin{bmatrix}\lambda t_0\mat{\zeta}{}{}{\zeta^{-1}}&\\
  &^tt_0^{-1}\mat{\zeta^{-1}}{}{}{\zeta}\,\end{bmatrix}\Big)(\zeta-\zeta^{-3})
  \lambda^{-4}\,d\zeta\,d\lambda,
\end{align}
where $t_0 \in \GL(2,\R)^+$ is such that $T^1(\R) = t_0 \SO(2) t_0^{-1}$.
Let us first consider the case $D \equiv 0$ mod $4$, in which case
$S(-D) = \mat{\frac D4}{0}{0}{1}$ and
$t_0 = \mat{2^{\frac 12} D^{-\frac 14}}{}{}{2^{-\frac 12} D^{\frac 14}}$.
From Sect.\ 4.4 of \cite{P-S1}, we have for $\lambda > 0$
\begin{equation}\label{whittakerfirstformulaeq}
 W_\infty\Big(\eta \begin{bmatrix}\lambda t_0\mat{\zeta}{}{}{\zeta^{-1}}&\\
  &t_0^{-1}\mat{\zeta^{-1}}{}{}{\zeta}\end{bmatrix},s\Big)  = c(1) \Big|\lambda D^{-\frac12}
  \big(\frac{\zeta^2+\zeta^{-2}}{2}\big)^{-1}\Big|^{3(s+\frac12)} W_{\frac l2,
  \frac{ir}2}\big(4 \pi \lambda D^{1/2}\frac{\zeta^2+\zeta^{-2}}2\big),
\end{equation}
where $c(1)$ is the first Fourier coefficient of $f_l$.
Using the results obtained in Sect.\ \ref{finding-good-vector-sec},
we will state the formula for the Bessel function appearing in (\ref{archintegral1eq}).
Observe that $\Lambda_\infty \equiv 1$. By arguments from Sect.\
\ref{besselmodelssec} on change of Bessel models, and the fact
that $\pi_{\bf F}$ has trivial central character, we obtain
\begin{equation}\label{Bessel-formula}
 B_{\bar{\phi}_\infty}\Big(\begin{bmatrix}\lambda t_0\mat{\zeta}{}{}{\zeta^{-1}}&\\
  &^tt_0^{-1}\mat{\zeta^{-1}}{}{}{\zeta}\,\end{bmatrix}\Big) =
  B_{l,l}(h(\lambda D^{\frac 12}2^{-1}, \zeta,0,0)),
\end{equation}
where $B_{l,l}$ is the vector in the $1$-dimensional $K$-type
obtained in Section \ref{finding-good-vector-sec}. We
computed the precise formula for $B_{l,l}$ for $\dim(V) = n =
3,5,7,9$. In each case, the formula has the shape
\begin{equation}\label{good-vector-general-formula}
 B_{l,l}(\lambda D^{\frac 12}2^{-1}, \zeta) =
  \Big(\sum\limits_{j=0}^{\big[\frac{n-1}4\big]}\sum\limits_{k=2j}^{\frac{n-1}2} c_{k,j}
  \Big(\frac{\lambda D^{\frac12}}2\Big)^{l-k}x^{k-2j}\Big)
   e^{-4 \pi \lambda\frac{D^{1/2}} x},
\end{equation}
where $x=(\zeta^2+\zeta^{-2})/2$ and $c_{k,j}$ are real constants
depending on $k, j, l, n$. Hence
$$
 Z_{\infty}(s) =\sum\limits_{j=0}^{\big[\frac{n-1}4\big]}
 \sum\limits_{k=2j}^{\frac{n-1}2}c_{k,j}\,Z_{\infty}^{k,j}(s),
$$
where
\begin{align*}
 Z_{\infty}^{k,j}(s) &=\pi \Big(\frac{D^{\frac12}}2\Big)^{l-k}c(1)
 \int\limits_0^{\infty} \int\limits_1^{\infty}
  \lambda^{3(s+\frac 12)+l-k}D^{-\frac 32(s+\frac 12)} x^{-3(s+\frac12)+k-2j}
  W_{\frac l2, \frac{ir}2}(4\pi\lambda D^{\frac12}x)e^{-2\pi\lambda D^{\frac 12}x}
  \lambda^{-4}\,dx\,d\lambda\\
 &= \frac{c(1) 2^{-6s+3-3l+3k} D^{-3s} \pi^{-3s-l+k+\frac
 52}}{6s+l-2k-1+2j}\frac{\Gamma(3s+l-k-1+\frac{ir}2)
 \Gamma(3s+l-k-1-\frac{ir}2)}{\Gamma(3s+\frac l2 - k - \frac
 12)}\\
 &=Q_{k,j}(s) c(1) 2^{-6s+3-3l} D^{-3s} \pi^{-3s-l+\frac 52}
  \frac{\Gamma(3s+l-1+\frac{ir}2)\Gamma(3s+l-1-\frac{ir}2)}{\Gamma(3s+\frac l2 - \frac 12)}
\end{align*}
with
$$
 Q_{k,j}(s) = \frac{2^{3k} \pi^{k}}{6s+l-2k-1+2j}
 \prod\limits_{t=1}^k\frac{3s+\frac l2- \frac 12 -t}
 {(3s+l-t-1+\frac{ir}2)(3s+l-t-1-\frac{ir}2)}.
$$
Thus, we get
\begin{equation}\label{arch-int-final-exp}
 Z_{\infty}(s) = \Big(\sum\limits_{j=0}^{\big[\frac{n-1}4\big]}
 \sum\limits_{k=2j}^{\frac{n-1}2} c_{k,j} Q_{k,j}(s)\Big) c(1)
 2^{-6s+3-3l} D^{-3s} \pi^{-3s-l+\frac 52}\frac{\Gamma(3s+l-1+\frac{ir}2)
  \Gamma(3s+l-1-\frac{ir}2)}{\Gamma(3s+\frac l2 - \frac 12)}.
\end{equation}
Recall that we did the above calculations under the assumption
that $D \equiv 0 \pmod{4}$. Following \cite{Fu} or the methods
from Sect.\ 4.4 of \cite{P-S1}, we get the same formula for
$Z_{\infty}(s)$ for $D \equiv 3 \pmod{4}$. Using the
non-archimedean calculations from \cite{Fu}, \cite{P-S1} and
\cite{P-S2} we now get the following global theorem.
\begin{theorem}\label{main-global-thm}
 Let $\rho = \det^l \otimes \rho_0$, where $l$ is an even
 positive integer and $\rho_0$ is the $n$-dimensional, irreducible representation
 of $\PGL(2,\C)$ with $n = 3,5,7$ or $9$. Let ${\bf F} \in S_\rho(\Gamma_2)$ be a cuspidal
 vector valued Siegel eigenform of degree $2$ satisfying the two assumptions
 in Sect.\ \ref{vvsec}. Let $L=\Q(\sqrt{-D})$, where $D$ is as in Assumption 1.
 Let $N = \prod p^{n_p}$ be a  positive integer. Let $f$ be a Maa{\ss}
 Hecke eigenform of weight $l_1 \in \Z$ with respect to $\Gamma_0(N)$. If $f$ lies in a
 holomorphic discrete series with lowest weight $l_2$, then assume that $l_2 \leq l$.
 Then the integral (\ref{global integral calculation}) is given by
 \begin{equation}\label{integral-l-fn-formula}
  Z(s,\Lambda) = \kappa_N Z_{\infty}(s) \frac{L(3s+\frac 12, \pi_{\bf F} \times
  \tau_f)}{\zeta(6s+1) L(3s+1, \tau_f\times \AI(\Lambda))},
 \end{equation}
 where
 \begin{align*}
  \kappa_N &=\prod\limits_{p | N} \frac{p-1}{p^{3(n_p-1)}(p+1)(p^4-1)}
    (1-\Big(\frac L{p}\Big)p^{-1})p^{n_p}   (1-p^{-6s-1})^{-1}
   \prod\limits_{p^2 | N} L_p(3s+1, \tau_p \times\AI(\Lambda_p)),\\
  \Big(\frac Lp\Big)&=\left\{\begin{array}{l@{\qquad\text{if }p}l}
   -1&\text{ is inert in }L,\\
   0&\text{ ramifies in }L,\\
   1&\text{ splits in }L.\end{array}\right.
\end{align*}
$Z_{\infty}(s)$ is as in (\ref{arch-int-final-exp}) and
$\AI(\Lambda)$ is the representation of $\GL(2,\A)$ obtained by
automorphic induction from $\Lambda$.
\end{theorem}

\end{document}